\documentclass[11pt,reqno]{amsart}
\usepackage[margin=1in]{geometry}
\makeatletter
\def\section{\@startsection{section}{1}%
	\z@{.7\linespacing\@plus\linespacing}{.5\linespacing}%
	{\bfseries
		\centering
}}
\def\@secnumfont{\bfseries}
\makeatother

\usepackage{amsmath,amssymb,amsthm,graphicx,amsxtra, setspace}
\usepackage[utf8]{inputenc}
\usepackage{charter}
\usepackage{mathrsfs}
\usepackage{alltt}
\usepackage{stackengine}
\usepackage{relsize}
\usepackage{hyperref}
\usepackage{aliascnt}
\usepackage{tikz}
\usepackage{mathtools}
\usepackage{multicol}
\usepackage{upgreek}
\usepackage{graphicx,type1cm,eso-pic,color}
\allowdisplaybreaks
\usepackage{scalerel,stackengine}
\stackMath
\usepackage{setspace}
\newtheorem{theorem}{Theorem}[section]
\newtheorem*{theorem*}{Theorem}

\newaliascnt{lemma}{theorem}
\newtheorem{lemma}[lemma]{Lemma}
\aliascntresetthe{lemma} 

\newaliascnt{proposition}{theorem}
\newtheorem{proposition}[proposition]{Proposition}
\aliascntresetthe{proposition}

\newaliascnt{assumption}{theorem}

\aliascntresetthe{assumption}

\newaliascnt{auxiliary}{theorem}

\aliascntresetthe{auxiliary}

\newaliascnt{corollary}{theorem}
\newtheorem{corollary}[corollary]{Corollary}
\aliascntresetthe{corollary}

\newaliascnt{definition}{theorem}
\newtheorem{definition}[definition]{Definition}
\aliascntresetthe{definition}

\newaliascnt{example}{theorem}

\aliascntresetthe{example}

\newaliascnt{remark}{theorem}
\newtheorem{remark}[remark]{Remark}
\aliascntresetthe{remark}

\newaliascnt{hypothesis}{theorem}

\aliascntresetthe{hypothesis}

\newaliascnt{property}{theorem}

\aliascntresetthe{property}

\usepackage[utf8]{inputenc}
\usepackage{amsmath}
\usepackage{enumerate}
\usepackage{bm}
\usepackage[T1]{fontenc}
\usepackage{amsthm}
\usepackage{amsfonts}
\usepackage{amssymb}
\usepackage{graphicx}
\usepackage{enumerate}
\usepackage[english]{babel}
\usepackage{comment}
\usepackage{amssymb}
\usepackage{amsthm}
\usepackage{xcolor}
\usepackage{dutchcal}
\usepackage{esint}

\usepackage[title]{appendix}
\usepackage{xcolor}

\usepackage{amssymb,amsmath,amsthm}
\usepackage{hyperref}
\usepackage{xcolor}
\usepackage{graphicx}
\newcommand{\RR}{\mathbb{R}}
\newcommand{\Om} {\Omega}

\newcommand{\ullb} {{\underline{u}}}

\newcommand{\ulwt} {{\widetilde{u}_{\lambda}}}

\newcommand{\Ellb} {{\underline{E}}}

\def\dis{d}
\def\loc{\text{loc}}

\def\Div{div}
\def\Osc{osc}

\newcommand{\al} {{\underline{\alpha}}}

\def\kg{K_{\gamma}}
\def\kge{K_{\gamma,\varepsilon}}
\def\fp{{-\Delta_{p}}}

\def\fps{(-\Delta)_{p}^{s}}
\def\fqs{(-\Delta)_{q}^{s}}

\def\va{\varepsilon}
\def\tu{\widetilde{u}}

\def\th{\widetilde{h}}
\def\tf{\widetilde{f}}
\def\tw{\widetilde{w}}
\newcommand{\wop}{W^{1,p}_0}
\def\pv{P.V.}

\newcommand{\wsqq}{W^{s,q}}
\newcommand{\wsq}{W^{s,q}_0}

\newcommand{\wopo}{W^{1,p}_0(\Omega)}
\newcommand{\wsqo}{W^{s,q}_0(\Omega)}
\newcommand{\uep}{u_{\varepsilon}}
\newcommand{\wro}{ \overline{\omega}_{\rho}}
\newcommand{\T}{\mathcal{T}}
\newcommand{\Sl}{\mathcal{S}}
\newcommand{\w}{\mathcal{W}}
\newcommand{\wo}{\mathcal{W}_0}
\newcommand{\wpp}{W^{1,p}}
\newcommand{\ck}{\mathfrak{c}}
\newcommand{\xo}{x_{0}}

\newcommand{\wu}{\widetilde{u}}
\newcommand{\wf}{\widetilde{f}}

\newcommand{\wO}{\widetilde{\Omega}}
\def\supp{supp}

\newcommand{\wkg}{\widetilde{K}_{\gamma}}
\def\uo{u_{0}}
\def\wuo{\widetilde{u}_{0}}

\usepackage[title]{appendix}
\usepackage{scalerel}[2014/03/10]
\usepackage{stackengine}
\def\dashint{\,\ThisStyle{\ensurestackMath{%
  \stackinset{c}{.2\LMpt}{c}{.5\LMpt}{\SavedStyle-}{\SavedStyle\phantom{\int}}}%
  \setbox0=\hbox{$\SavedStyle\int\,$}\kern-\wd0}\int}
\def\intavg{\,\ThisStyle{\ensurestackMath{%
  \stackinset{c}{.2\LMpt}{c}{.2\LMpt}{\SavedStyle-}{\SavedStyle\phantom{\int}}}%
  \setbox0=\hbox{$\SavedStyle\int\,$}\kern-\wd0}\int}

\DeclareMathOperator*{\essinf}{ess\,inf}

\newenvironment{sketch}{%
  \proof}{\endproof}
\newcommand{\Addresses}{{
		\footnote{
				\footnotesize
\noindent \textsuperscript{1}Indian Institute of Science Education and Research, Thiruvananthapuram 695551, India  \par\nopagebreak 
   \noindent \textsuperscript{2}Universite de Pau et des Pays de l’Adour, LMAP (UMR E2S-UPPA CNRS 5142) Bat. IPRA, Avenue de l’Universite F-64013 Pau, France  \par\nopagebreak 
\noindent \textsuperscript{A}\textit{e-mail:} \texttt{dhanya.tr@iisertvm.ac.in}.
\noindent \textsuperscript{B}\textit{e-mail:} \texttt{ritabrata20@iisertvm.ac.in}.
\noindent \textsuperscript{C}\textit{e-mail:} \texttt{jacques.giacomoni@univ-pau.fr}.

			 \noindent \textsuperscript{*}Corresponding author.

			\medskip\noindent
			{\bf Acknowledgments:} 
  R. Dhanya was supported by SERB MATRICS grant MTR/2022/000780 when this work was being carried out. J. Giacomoni was partially funded by IFCAM (Indo-French Centre for Applied Mathematics), IRL CNRS 3494. Ritabrata Jana expresses gratitude for the financial assistance provided by the Prime Minister's Research Fellowship during the execution of this research.  Furthermore, Ritabrata Jana extends heartfelt thanks to the University of Pau for for its hospitality and support during the conduct of a significant portion of this work.
			
}}}
\begin{document}
\title[Regularity of Mixed Local Nonlocal Problem with Singular Data ]{ Interior and Boundary Regularity of Mixed Local Nonlocal Problem with Singular Data and Its Applications	\Addresses	}
	\author[ R. Dhanya ]
	{ R. Dhanya\textsuperscript{1,A}} 
    \author[Jacques Giacomoni]
	{Jacques Giacomoni\textsuperscript{2,A,*}}
 	\author[Ritabrata Jana]
	{Ritabrata Jana\textsuperscript{1,B}}

\maketitle
\begin{abstract}
In this article, we examine the H\"older regularity of solutions to equations involving a mixed local-nonlocal nonlinear nonhomogeneous operator $\fp + \fqs$ with singular data, under the minimal assumption that $p> sq$. The regularity result is twofold: we establish interior gradient H\"older regularity for locally bounded data and boundary regularity for singular data. We prove both boundary H\"older and boundary gradient H\"older  regularity depending on the degree of singularity. Additionally, we establish a strong comparison principle for this class of problems, which holds independent significance. As the applications of these qualitative results, we further study sublinear and subcritical perturbations of singular nonlinearity.
\end{abstract}
\\
\section{Introduction}
In this article, we focus on the nonhomogeneous mixed local-nonlocal problem involving a singular term:
\begin{equation}\label{maineqn}
    \begin{aligned}
         \fp u +  \fqs u &= \frac{\kg(x)}{u^{\delta}}+F(x,u), \quad u>0 \text{ in } \Omega,
        \\
        u &=0 \text{ in } \RR^N \setminus \Omega
    \end{aligned}
\end{equation}
where  $p,q> 1,$ $p\geq sq$ and $\gamma,\delta>0.$ Here $\kg$ is a nonnegative weight function  which satisfies the  condition:
  \begin{equation*}
 \kg(x)\in L^{\infty}_{\loc}(\Omega) \;\;\; \text{ and }   \;\;\;\; C_1 \dis_{\Omega}^{-\gamma}(x)\leq \kg(x) \leq C_2 \dis_{\Omega}^{-\gamma}(x) \text{ in } \Omega
\end{equation*}
for some $C_1, C_2>0.$ Here $F$ is a Carath\'eodory function and $\Omega \subset \mathbb{R}^N$ is an open, bounded domain with a smooth boundary. We define the distance function $\dis: \RR^N \rightarrow \RR_{+}$ as 
\begin{equation*}
    \dis(x):=\dis_{\Omega}(x):=\underset{y \in \Om^c}{\inf \thinspace } |x-y|.
\end{equation*} The local part of the operator is defined as $ \fp u:= -\nabla.(|\nabla u|^{p-2}\nabla u).$ For the nonlocal part, the fractional $q$-Laplacian, denoted by $\fqs$, is defined as 
\begin{equation*}\label{o11}
(-\Delta)_{q}^su(x)= 2\thinspace PV. \int_{\mathbb{R}^N}\frac{|u(x)-u(y)|^{q-2}(u(x)-u(y)}{|x-y|^{N+sq}} \,dy \qquad x\in \Omega,
\end{equation*}
up to a suitable normalization constant. 
   The study of mixed local-nonlocal operators has a long history, dating back several decades \cite{BCP68, Can86}. 
   Operators of this type naturally emerge from the interplay of two stochastic processes operating at different scales, specifically, a classical random walk and a Lévy flight. In essence, when a particle transitions between these two processes based on a defined probability distribution, the corresponding limit diffusion equation is governed by the aforementioned operator.  For a detailed discussion of this phenomenon, please see the appendix in \cite{DPV23}. These operators also play an essential role in various applications, including the biological sciences, as highlighted in \cite{DV21} and related references, as well as in the study of heat transport in magnetized plasmas \cite{BD13}, among others. In this manuscript, we explore the interior and the boundary regularity results for mixed nonhomogeneous operators involving singular data. Towards the end of the manuscript we also establish a comparison principle and obtain other existence and uniqueness results. 
\par

To lay the groundwork for the paper, we begin by examining the regularity results for equations of the form \(\mathcal{L}u = f(x)\) in \(\Omega\) with \(u = 0\) in \(\Omega^c\). For local operators such as the \(p\)-Laplacian, it is well established that if \(f\) is a bounded measurable function, the solution \(u\) belongs to \(C^{1,\alpha}(\overline{\Omega})\) see \cite{LU68,Lib88}. However, in the nonlocal setting, where \(\mathcal{L}\) represents operators like the fractional Laplacian, fractional \(p\)-Laplacian, or similar nonlocal operators, the maximal regularity achievable up to the boundary is H\"older continuity (see \cite{RS14,IMS16,GKS23}). More specifically, authors have obtained a uniform \(C^\alpha\) bound for solutions \(u\) in terms of the \(L^\infty\) bound of the given data \(f\). Moreover, an improved regularity result is obtained in Iannizzotto et al.\cite{IMS20} and Ros-Oton et al.\cite{RS14} where the authors proved that the solution $u$ belongs to a weighted H\"older space. In contrast, for mixed local-nonlocal operators, $\mathcal{L}u=$ \(-\Delta_p u + (-\Delta)^s_q u = f(x)\), when \(p \geq sq\), the situation is significantly different. Under appropriate integrability conditions on \(f\) and suitable boundary conditions, using D. Filippis et al.\cite{FM22} and  Antonini et al.\cite{AC23} it can be seen that the solution \(u\)  belongs to \(C^{1,\alpha}(\mathbb{R}^n)\). In this article, our primary focus is to establish similar results under the assumptions that \(f\) is locally bounded in \(\Omega\) and exhibits a behavior of \(d(x)^{-\gamma}\) as \(x\) approaches the boundary of \(\Omega\) subject to homogeneous boundary condition. Once the interior and boundary regularity results for \(-\Delta_p u + (-\Delta)^s_q u = f(x)\) are established, we shift our focus to nonlinear problems involving singular terms of the prototype $ -\Delta_p u + (-\Delta)^s_q u = u^{-\delta}$ for $\delta>0.$

The study of singular problems for the Laplace operator dates back to the 1960's \cite{FM60}. However, the seminal work of Crandall, Rabinowitz, and Tartar \cite{CRT77} drew significant attention to this class of singular elliptic PDEs. Interest in the topic grew further following the renowned paper of Lazer and McKenna \cite{LM91}, where it was demonstrated that for strongly singular problems, classical solutions may not lie within the natural Sobolev space. Later, these results were also studied and generalized for nonlinear operators such as $\fp$ and $\fp -\Delta_q$ with singular terms, see \cite{CST16,PS06,LP23}. For these operators, the interior regularity result can be inferred from the classical text of Ladyenskaja-Ural'ceva \cite{LU68}, but the boundary regularity has to be worked out. Giacomoni et al. \cite{GST07, GKS21} extensively studied the Sobolev and H\"older regularity of these operators involving singular terms near the boundary. Once the regularity results are established, the correct functional framework is set to analyze several existence/nonexistence, uniqueness and multiplicity results for nonlinear problems using topological or variational techniques
(see \cite{DIJ23,DI23,Aro22,LP23}). Also, several qualitative properties of the solution such as comparison principle \cite{GST07}, Picone's identity\cite{BT20}, Sobolev versus H\"older minimizers results \cite{GS10} etc can also be deduced with the help of $C^{1,\alpha}$ boundary regularity result.

\par
Motivated by the study of singular elliptic problems for local operators, several researchers \cite{BDMP15, CMSS17, CMST23, AGS18} have extended their investigations to similar problems involving nonlocal operators such as $(-\Delta)^s$ and $\fps.$ Primary aim of these works were to establish the existence or multiplicity of solutions within suitable Sobolev spaces using variational methods. While interior regularity results for these problems can often be inferred from \cite{BLS18}, it was first in \cite{AGW21}  that the H\"older-type boundary regularity for singular problems was  addressed. Later, Giacomoni et. al.\cite{GKS23} established the $C^{\alpha}(\overline{\Omega})$ regularity results for fractional $(p,q)$ Laplace operator with singular nonlinearity. One of the major focuses of this article is to establish the existence and H\"older regularity of solutions of mixed local-nonlocal equations involving a singular nonlinearity.
\par
We shall now provide a brief overview of the mixed local nonlocal boundary value problems of the type $\fp u+(-\Delta)^s_q u= u^{-\delta}$ in $\Omega$ satisfying $u=0$ in $\Omega^c.$ In the specific case of a linear mixed local-nonlocal operator with a singular term (i.e., $p = q = 2$ and $\delta > 0$), Sobolev regularity and boundary behavior were obtained in \cite{AR23} using the Green's function approach. In \cite{GU22}, existence, uniqueness, and Sobolev regularity of weak solutions of homogeneous case i.e. when $p = q$ are discussed. Through approximation methods and local regularity results, existence and regularity were explored for both constant and variable singular exponents in \cite{GKK24}.  Furthermore, multiplicity results for subcritical perturbations, particularly in the case $p = q = 2$ in \eqref{maineqn}, were discussed in \cite{Gar23} and \cite{BD24a}.

\par
Building upon the prior findings, in this article, we first establish local boundedness for solutions with locally bounded data under the condition $q>p.$ This partially addresses a broader and more challenging problem as outlined in \cite[Section 8.3]{FM22}. We begin with Caccioppoli-type estimates, which cannot be directly borrowed from existing results due to the non-homogeneity of the operator and the lack of embedding in the Sobolev space. From these estimates, we derive the desired local boundedness result, and by applying \cite[Theorem 6.1]{FM22}, we obtain interior $C^{1,\alpha}_{\loc}$ regularity.

 A further novelty of this work lies in adapting the perturbative methods developed for mixed local-nonlocal operators with bounded data in \cite{FM22} and \cite{AC23} to the setting of singular data. Employing the barrier method, which is widely known for regularity results in nonlocal operators, we first establish the H\"older regularity of weak solutions up to the boundary. Leveraging the boundary behavior, we modify the perturbative method described in \cite{AC23} to obtain the gradient H\"older boundary regularity result. 
 \par 
Next, we also analyze the regularity of solutions to the mixed local-nonlocal problem involving a singular nonlinearity, specifically \eqref{maineqn}, where \(F(x, u) = f(x)\) with \(f(x) \in L^{\infty}(\Omega)\). Our study focuses on a detailed examination of the behavior of barrier functions when acted upon by mixed local-nonlocal operator. The analysis begins with establishing the existence and uniqueness of the weak solution to an auxiliary approximated problem, as formulated in \eqref{auxpb}. By utilizing Hardy inequalities in both local and fractional cases, we derive the Sobolev regularity of the solutions. We identify the appropriate exponents by considering  distinct cases for the parameters $\gamma,\delta,s$ and $q.$ We then investigate the boundary behavior of the solutions through the use of barrier functions and the property of local boundedness, thereby establishing H\"older regularity.

\par
Once the regularity results are obtained we focus on the application part. We establish the strong comparison principle (SCP) for mixed local-nonlocal operators with singular nonlinearity ($\gamma=0,$ $0<\delta<1$ in \eqref{maineqn}), a result that holds independent significance. Once gradient H\"older regularity is established for singular nonlinearities, it is natural to assume that the weak solution lies within the gradient H\"older space. Leveraging estimates from \cite{Jar18}, we prove the SCP for nonlinear mixed local-nonlocal operators with singular terms. 

We conclude this article by presenting two applications of the qualitative properties established earlier. Namely, we examine the case 
$F(x,u) = u^l$ in \eqref{maineqn} under two conditions: one with sublinear perturbations and another with superlinear but subcritical perturbations. A major challenge in both cases is proving uniform boundedness for the nonlinear problem and  detailed statements of these applications are provided in Theorem \ref{extncuniqsub}, and Theorem \ref{extsubcr}.

We shall now present the main novelty of this article.
\begin{enumerate}
 \setlength \itemsep{1em} 
    \item We establish regularity results for mixed local and nonlocal operators of the type $\fp+\fqs$ without imposing additional restrictions, other than assuming $p>sq.$ In the light of remark 3 of Mingione-De Filippis\cite{FM22} the case $p>q$ is comparatively easier to handle. This is due to embedding results for fractional Sobolev space where $W^{1,p}\cap W^{s,q}= W^{1,p}$ and hence the $p$ Laplace operator dominates the fractional part. In contrast, the focus of our analysis is on the more intricate case $(sq< p <q)$ where gradient norm cannot dominate the fractional norm,  introducing significant analytical challenges.

    \item Until recently, research on singular elliptic problems involving $\fp+\fqs$ was only focusing on Sobolev type regularity results. This work represents one of the first attempts to establish a  H\"older type regularity results.   The results we establish for local solutions are new, even for some of the homogeneous mixed operators, i.e., when $p = q$.
    \item For the reasons mentioned in point (1), the regularity result we have obtained is for any $p>sq,$ which is the best range we can expect, as seen in  Mingione-De Filippis\cite{FM22} or Antonini-Cozzi\cite{AC23}. The results we establish play a fundamental role in understanding nonlinear elliptic problems and proving existence and multiplicity results. Furthermore, this work explores applications of the derived regularity, which includes a comparison principle and the existence of solutions for certain nonlinear elliptic problems.
\end{enumerate}
Towards the completion of this article we came across a recent preprint of K.Bal and S.Das \cite{BD24b} where they consider similar problems but with restricted conditions such as $ q\leq p.$ However, their definition of a weak solution assumes weaker tail conditions compared to ours. In the mixed local-nonlocal case, when $q\leq p$, the local operator becomes dominant and in this framework, their result provides sharp boundary estimates. On the other hand, when $sq<p<q$, both operators contribute significantly. In this setting, under the presence of singular data, our result is the first of its kind.

\par
The article is structured as follows: Section \ref{prlm} introduces essential definitions and presents the main results. In Section \ref{regsingular}, we discuss interior regularity for local data and establish the H\"older and gradient H\"older regularity of solutions up to the boundary. Section \ref{regsingnonlinear} addresses Sobolev and H\"older regularity for singular nonlinearities. In Section \ref{comp}, we prove comparison principles for solutions. Applications of the main results are explored in Sections \ref{picone} and \ref{subcrtical}: Section \ref{picone} covers sublinear perturbations, while Section \ref{subcrtical} examines subcritical, superlinear perturbations.
\par
\textbf{Notations:}
Unless stated otherwise, $k,M, C$ etc. represent generic positive constants whose meaning can be different in  the even same line. Throughout this article we have assumed $p,q>1$ and $p> sq.$ Given $a\in \RR,$ we denote $a_{+}:=\max\{a,0\}.$ We use the notation $[a]^{t}:=|a|^{t-1}a$ for any $a\in \RR, t>0.$ Let $B_R(x_0)$ denote the open ball with radius $R>0$ and centre at $x_0\in \mathbb{R}^N.$

We omit denoting the centre when it is not necessary. \\
For $S\subset\RR^{2N},$ we denote
\begin{equation}
    \begin{aligned}
        A_{t}(u,v,S):=\int_{S}\frac{|u(x)-u(y)|^{t-2}(u(x)-u(y))(v(x)-v(y))}{|x-y|^{N+sq}} \thinspace dx\thinspace dy.
    \end{aligned}
\end{equation}
We say $f\lesssim g$  in $\Omega$ if there exists a $C(N,s,p,q,\Omega)>0$ such that $f\leq C g$ in $\Omega.$ Other quantities $p,q$ and $\kg$ remain fixed as described in the introduction.

\section{Preliminaries And Main Results}\label{prlm}

\subsection{Function Space}
We recall that for $E \subset \mathbb{R}^N$, the Lebesgue space $L^t(E)$, $1 \leq t < \infty$, is defined as the space of $t$-integrable functions $u: E \to \mathbb{R}$ with the finite norm 
\begin{equation*}
    \|u\|_{L^t(E)} = \left( \int_E |u(x)|^t \, dx \right)^{1/t}.
\end{equation*}
 The Sobolev space $W^{1,t}(\Omega)$, for $1 \leq t < \infty$, is defined as the Banach space of locally integrable weakly differentiable functions $u: \Omega \to \mathbb{R}$ equipped with the following norm 
\begin{equation*}
    \|u\|_{W^{1,t}(\Omega)} = \|u\|_{L^t(\Omega)} + \|\nabla u\|_{L^t(\Omega)}.
\end{equation*}
The space $W^{1,t}_0(\Omega)$ is defined as the closure of the space $C_c^\infty(\Omega)$ in the norm of the Sobolev space $W^{1,t}(\Omega)$, where $C_c^\infty(\Omega)$ is the set of all smooth functions whose supports are compactly contained in $\Omega$. For a measurable function $u:\mathbb{R}^{N}\rightarrow \mathbb{R}$, we define Gagliardo seminorm
\begin{equation*}
	[u]_{s,t}:=[u]_{W^{s,t}(\mathbb{R}^N)}:=  \left(\int_{\mathbb{R}^N \times \mathbb{R}^N} \dfrac{|u(x)-u(y)|^{t}}{|x-y|^{N+st}}\,dx\, dy\right)^{1/t}
\end{equation*}
for $1<t<\infty$ and $0<s<1.$
We consider the space $W^{s,t}(\mathbb{R}^N)$ defined as  
\begin{equation*}
	W^{s,t}(\mathbb{R}^{N}):= \left \{u \in {L}^{t}(\mathbb{R}^{N}):[u]_{s,t}<\infty  \right\}.
\end{equation*}
The space $W^{s,t}(\mathbb{R}^{N})$ is a Banach space with respect to the norm 
\begin{equation*}
\|u\|_{{W^{s,t}(\mathbb{R}^{N})}}=\left( \|u\|^{t}_{L^{t}(\mathbb{R}^{N})} + [u]^{t}_{{W^{s,t}(\mathbb{R}^{N})}}\right)^{\frac{1}{t}} . 	
\end{equation*}
A comprehensive examination of the fractional Sobolev Space and its properties are presented in \cite{DPV12}. To address the Dirichlet boundary condition, we naturally consider the space $W^{s,t}_0(\Omega)$ defined as
\begin{equation*}
	W_{0}^{s,t}(\Omega):= \left \{u \in W^{s,t}(\mathbb{R}^{N}):u=0\medspace\text{in}\medspace \mathbb{R}^{N}\setminus \Omega \right\}.
\end{equation*}
This is a separable, uniformly convex Banach space endowed with the norm $\|u\|=	\|u\|_{{W^{s,t}(\mathbb{R}^{N})}}$. 
Now, we define the local spaces as 
\begin{equation*}
    \begin{aligned}
        &L^t_{\text{loc}}(E)=\left\{u:\Omega\rightarrow\RR: u \in L^t(K) \text{ for every } K \subset\subset E\right\},
        \\
        &W^{1,t}_{\text{loc}}(\Omega) = \left\{ u : \Omega \to \mathbb{R} : u \in L^t(K), \int_K |\nabla u|^t \, dx < \infty, \, \text{for every } K \subset\subset \Omega \right\},
        \\
        &W^{s,t}_{\text{loc}}(\Omega) = \left\{ u : \Omega \to \mathbb{R} : u \in L^t(K), \int_K \int_K \frac{|u(x) - u(y)|^t}{|x - y|^{N + st}} \, dx \, dy < \infty, \, \text{for every } K \subset\subset \Omega \right\}.
    \end{aligned}
\end{equation*}
Now, for $N > t$, we define the critical Sobolev exponent as $t^{*} = \frac{Nt}{N - t},$
and we get the following embedding result for any bounded open subset $\Omega$ of class $C^1$ in $\mathbb{R}^N:$ there exists $C \equiv C(N, \Omega) > 0$ such that for all $u \in C_c^{\infty}(\Omega)$,  
\begin{equation*}
    \begin{aligned}
        \|u\|_{L^{t^*}(\Omega)} \leq C \int_{\Omega} |\nabla u|^t \, dx.
    \end{aligned}
\end{equation*}
Moreover, the inclusion map
\begin{equation*}
    \begin{aligned}
        W_0^{1,t}(\Omega) \hookrightarrow L^r(\Omega)
    \end{aligned}
\end{equation*}
is continuous for $1 \leq r \leq t^*$, and the above embedding is compact except for $r = t^*$. Also the embedding $W^{s,t}_{0}(\Omega)\hookrightarrow L^{r}(\Omega)$ is continuous for $1\leq r\leq t^{*}_{s}:=\frac{Nt}{N-ts}$ and
compact for  $1\leq r < t^{*}_{s}$. 
 Due to continuous embedding of $W^{s,t}_{0}(\Omega)\hookrightarrow L^{r}(\Omega)$ for $1\leq r\leq t^{*}_{s}$, we define the equivalent norm on $W^{s,t}_{0}(\Omega)$ as
\begin{equation*}
	\|u\|_{W^{s,t}_{0}}:=\left(\int_{\mathbb{R}^N \times \mathbb{R}^N} \dfrac{|u(x)-u(y)|^{t}}{|x-y|^{N+st}}\,dx\, dy\right)^{1/t}.
\end{equation*} 
The dual space of $W^{s,t}_{0}(\Omega)$ is denoted by $W^{-s,t'}(\Om)$ for $1<t<\infty.$
\\ 
Throughout this paper, unless stated otherwise, no specific relationship between  $p$  and $q$ is assumed beyond  $p > sq$. Consequently, there may be a lack of continuous embedding between $\wop$  and  $\wsq$. To address this, when analysing the weak solution associated with the operator $ \fp + \fqs$, we consider
\begin{equation*}
    \begin{aligned}
        \mathcal{W}(\Omega) = W^{1,p}(\Omega) \cap W^{s,q}(\Omega),
    \end{aligned}
\end{equation*}
with the norm  $\|\cdot\|_{\mathcal{W}(\Omega)} = \|\cdot\|_{W^{1,p}(\Omega)} + \|\cdot\|_{W^{s,q}(\Omega)}$ . To incorporate the zero Dirichlet boundary condition, we also define 
\begin{equation*}
    \begin{aligned}
       \mathcal{W}_0(\Omega) = W_0^{1,p}(\Omega) \cap W_0^{s,q}(\Omega).
    \end{aligned}
\end{equation*}
The space $\mathcal{W}^{\prime}(\Om)$ is denoted for the dual of the $\w(\Om).$ 
We introduce the tail space as
\begin{equation*}
    \begin{aligned}
        L^{m}_{\alpha}(\RR^{N}):=\left\{u\in L^{m}_{\loc}(\RR^N): \int_{\RR^{N}} \frac{|u(x)|^{m}}{(1+|x|)^{N+\alpha}}\thinspace dx<\infty\right\}.
    \end{aligned}
\end{equation*}
Next we define the following space for defining the local solution:
\begin{equation*}
    \begin{aligned}
        \mathcal{X}(\Omega):= \wpp_{\loc}(\Omega)\cap\wsqq_{\loc}(\Omega)\cap L^{q}_{sq}(\RR^N).
    \end{aligned}
\end{equation*}
\subsection{Main Results}

We begin by defining a local solution, following the approach in \cite[Section 7]{FM22}:
\begin{definition}
    We will say $u\in\mathcal{X}(\Omega)$ be a local subsolution (supersolution) to the problem 
    \begin{equation}\label{localeqn}
        \begin{aligned}
            \fp u +\fqs u=f(x) \text{ in } \Omega
        \end{aligned}
    \end{equation}
    if for every non-negative $\varphi\in{ \mathcal{W}_0}(\Omega)$ with compact support the following holds: 
    \begin{equation*}
        \begin{aligned}
            \int_{\Omega} [ |Du|^{p-2}Du D\varphi -f \varphi ]+\int_{\RR^N}\int_{\RR^N} \frac{|u(x)-u(y)|^{q-2}(u(x)-u(y))(\varphi(x)-\varphi(y))}{|x-y|^{n+sq}} dx dy\leq(\geq) \ 0.
        \end{aligned}
    \end{equation*}
    We say $u$ is a solution if it is both a sub and supersolution.
\end{definition}

The following theorem provides an interior regularity result for solutions with local data:

\begin{theorem}
    \label{interiorregularity}
     Let $u\in\mathcal{X}(\Omega) $ be a local solution to \eqref{localeqn}. If $f\in L^{\infty}_{\loc}(\Omega)$ then 
     \begin{equation}\label{interiorholderreg}
         \begin{aligned}
             u \in C^{\alpha}_{\loc}(\Omega) \text{ for each } 0<\alpha<1.
         \end{aligned}
     \end{equation} 
     Moreover, there exists $\theta\in(0,1)$ such that $ u \in C^{1,\theta}_{\loc}(\Omega).$
\end{theorem}
We now examine the boundary behavior of weak solutions to the problem:
\begin{equation}\label{singregeqn1}
    \begin{aligned}
        \fp u +\fqs u&= f(x)\text{ in } \Omega,
        \\
        u&=0 \text{ in } \Omega^{c}
    \end{aligned}
\end{equation}
where $f\in L^{\infty}_{\loc}(\Omega)$ such that   \begin{equation*}
0<C_f \leq f(x) \leq C_{1}  \dis^{-\gamma}(x) \text{ in } \Omega.
\end{equation*}
\begin{theorem}\label{gradientholderregularity}
Let $0<\gamma<1$ and $u$ be a  weak solution to the problem
 \eqref{singregeqn1} belonging to $\wop(\Omega)\cap\wsq(\Omega)$. Then  $u \in C^{1,\theta}(\Bar{\Omega})$ for some $\theta\in (0,1).$
\end{theorem}
Now we discuss the regularity of the solution to the problem \eqref{maineqn}. We define the notion of a weak solution for problems involving singular data:
\begin{definition}\label{soldef}
    A function $u\in\mathcal{X}(\Omega) $ is said to be a sub-solution (super solution) of \eqref{maineqn} if the followings holds: 
    \begin{itemize}
        \item[i)] there exists $\theta\geq 1$ such that $u^{\theta}\in  \wopo\cap\wsqo,$
        \item[ii)] $\essinf_{K} u>0$ for all $K\subset\subset\Omega,$
        \item[iii)] for any $\varphi\in \wopo\cap\wsqo$ and $\varphi\geq 0,$
        \begin{equation*}
    \begin{aligned}
        \int_{\Omega} |\nabla u|^{p-2} \nabla u \nabla \varphi &+ \int_{\RR^{2N}}\frac{|u(x)-u(y)|^{q-2}(u(x)-u(y))(\varphi(x)-\varphi(y))}{|x-y|^{N+sq}} \thinspace dx\thinspace dy 
        \\
       & \leq (\geq) \int_\Omega (\frac{\kg(x)}{u^{\delta}}+F(x,u)) \varphi \;\;\;\;\;\mbox{  holds true.}
    \end{aligned}
\end{equation*}
    \end{itemize}
    We say that $u$ is a weak solution of \eqref{maineqn} if it is both a weak supersolution and subsolution to the problem.
\end{definition}
We now describe the boundary behavior of the solutions to the problem in two cases.
\begin{theorem}\label{gradhld}
     Let  $u$ be a weak solution to \eqref{maineqn} with $F(x,u)=f(x)\in L^{\infty}(\Omega)$ and $f(x)\geq 0$ where the exponents satisfy the conditions $\gamma+(\frac{1}{q}-s)(1-\delta)<1$ and $0<\gamma+\delta<1.$ Then there exists an $\alpha\in(0,1)$ such that $u\in C^{1,\alpha}(\Bar{\Omega}).$
\end{theorem}

\begin{theorem}\label{strongsingreg}
   Let $0<\gamma<p$ and $u$ be a minimal weak solution of \eqref{maineqn} with $F(x,u)=0$ where  $\gamma+\delta>1.$  Then there exists some $\alpha\in(0,1)$ such that $u\in C^{\alpha}(\Bar{\Omega}).$

\end{theorem}

\begin{remark}\label{cndgmdel}
    The condition $\gamma+(\frac{1}{q}-s)(1-\delta)<1$ in Theorem \ref{gradhld} can be omitted if $1\leq q \leq p$ or $q>\frac{1}{s}$ or $\gamma=0.$ 
\end{remark}

We now establish a form of the strong comparison principle for $C^{1,\alpha}$ solutions of mixed local-nonlocal problems.
\begin{theorem}[Strong Comparison Principle]\label{scp}
   Let us assume that $0<\delta<1, \lambda>0,$ $q>\frac{1}{1-s}$ and  $u,v\in C^{1,\alpha}(\Bar{\Omega})$ solves
\begin{equation}\label{uvsol}
    \begin{aligned}
        \fp u + \fqs u -\frac{\lambda}{u^{\delta}}&=f, \thinspace u>0 \text{ in } \Omega \text{ and } u|_{\Omega^c}=0,
        \\
         \fp v + \fqs v -\frac{\lambda}{v^{\delta}}&=g, \thinspace v>0 \text{ in } \Omega \text{ and } v|_{\Omega^c}=0
    \end{aligned}
\end{equation}
where $f,g$ are two continuous functions such that  $0\leq f<g$ pointwise. Then the following holds
\begin{equation*}
    \begin{aligned}
       0<u<v \text{ in } \Omega.
    \end{aligned}
\end{equation*}
\end{theorem}

\begin{remark}
    After establishing the strong comparison principle, a natural question arises regarding the validity of a Hopf-type lemma. Specifically, given \( u, v \in C^{1,\alpha}(\overline{\Omega}) \) that solve \eqref{uvsol}, one may ask whether the outward normal derivative satisfies  
\[
\frac{\partial}{\partial \nu}(v-u) < 0.
\]
This issue is addressed for local operator in \cite{GST07}. But for the nonlocal operator \(\fps\) a Hopf type lemma is established only for non-singular nonlinearities (see \cite{IMP23}).  The fractional counterpart of Hopf-type lemma for singular nonlinearities remains largely unexplored, especially for the fractional \( p \)-Laplacian \( (-\Delta)^{s}_{p} \) and even less is known for more complex operators involving mixed local-nonlocal cases. Notably, even for \( (-\Delta)^{s} \) with \( p = 2 \), questions of this nature are typically explored through the study of singular semipositone problems (see \cite{GMK19})—a framework that has yet to be extended to \( p \neq 2 \) in the nonlocal operator case. 
\end{remark}

As an application of these qualitative properties, we consider the case where $F(x,u) = u^l$ in \eqref{maineqn} or precisely we consider:
\begin{equation}\label{sublsuperineareqn}
    \begin{aligned}
        \fp u + \fqs u &=\frac{\lambda}{u^{\delta}}+u^{l} \text{ in }\Omega,
        \\
        u&=0 \text{ in } \Omega^c
        \end{aligned}
\end{equation}
 where $0<\delta<1$  and analyze two scenarios: one with sublinear perturbations i.e. $0<l< \min\{p-1,q-1\}$ and another with superlinear but subcritical perturbations i.e. $\max\{p-1,q-1\}<l<\min\{p^{*}-1,q_{s}^{*}-1\}$. For sublinear perturbations, we first establish a priori uniform bound for the solution. Alongside the established gradient Hölder boundary regularity, Picone's inequalities can be applied to mixed local-nonlocal operators, guaranteeing the uniqueness of the solution.

\begin{theorem}\label{extncuniqsub}
For any given $\lambda>0,$ there exists one and only one weak solution, $u$, to  \eqref{sublsuperineareqn} whenever $0<l< \min\{p-1,q-1\}.$ Furthermore, $u \in C^{1,\alpha}(\Bar{\Omega})$.
\end{theorem}
By standard minimization techniques applied to a suitable cut-off functional we prove the existence of a positive solution to the superlinear and sub-critical problem. We also prove a uniform $L^\infty$ bound for subcritical problems.

\begin{theorem}
   \label{extsubcr}
   Let $ \max\{p-1,q-1\}<l<\min\{p^{*}-1,q_{s}^{*}-1\}.$ Then there exists  $0<\Lambda\leq \infty$ such that for any $\lambda$ satisfying $0<\lambda<\Lambda,$ the equation $\eqref{sublsuperineareqn}$ admits a solution  $u \in C^{1,\alpha}(\Bar{\Omega})$ for some $\alpha \in (0,1).$  
\end{theorem}

\section{Regularity Results For Singular Data} \label{regsingular}
\subsection{Local Boundedness}
First we want to prove the local boundedness result for the local solution of \eqref{singregeqn1} when the provided data $f(x)$ belongs to $L^\infty_{\loc}(\Omega)$. Following the approach in \cite{BOS22}, we first establish the Caccioppoli inequality and then employ it to demonstrate local boundedness through the method of De Giorgi sequences. First we recall well known iteration lemma from \cite[Lemma 7.1]{Giu03}.
\begin{lemma}\label{itrlrm}
    Let $\{y_{i}\}_{i=0}^{\infty}$ be a sequence of nonnegative numbers satisfying $y_{i+1}\leq b_{1}b_{2}^{i}y_{i}^{1+\beta}$ for $i=0,1,2,\ldots$ for some constants $b_1,\beta>0$ and $b_{2}>1.$ If $y_{0}\leq b_{1}^{-1/\beta}b_{2}^{-1/\beta^{2}}$ then $y_{i}\rightarrow 0$ as $i\rightarrow \infty.$
\end{lemma}
Due to the work of \cite{FM22}, we only work on the case of $p\leq q.$ We begin by establishing the Caccioppoli estimates:
\begin{lemma}\label{cacciopoli}
    Let $1<p\leq q<\infty,$ $B_{2r}(x_0)\subset\subset\Omega$ be a ball  and $u$ be a local solution to \eqref{localeqn}. Define $w=(u-k)_{+}$ for $k>0.$ Then for any $\varphi\in C_{0}^{\infty}(B_{r})$ with $0\leq \varphi\leq 1$ we have
    \begin{equation*}
        \begin{aligned}
            \int_{B_{r}} \varphi^{q} |Dw|^{p} &+ \int_{B_{r}} \int_{B_{r}} \frac{|w(x)-w(y)|^{q}(\varphi^{q}(x)+\varphi^{q}(y))}{|x-y|^{N+sq}} 
            \\
            &\leq C\left[ \int_{B_{r}} w^{p} |D\varphi|^{p}+\int_{B_{r}} \int_{B_{r}} \frac{|\varphi(x)-\varphi(y)|^{q}(w(x)+w(y))^{q}}{|x-y|^{N+sq}} \right.
            \\& +\left. \left(\sup_{y\in \supp \varphi} \int_{\RR^{N}\setminus B_{r}} \frac{w(x)^{q-1}}{|x-y|^{N+sq}} \right) \int_{B_{r}} w \varphi^{q} +\int_{B_{2r}}|f|w\varphi^{q} \right]
        \end{aligned}
    \end{equation*}
    for some $C(N,s,p,q)>0.$
\end{lemma}
\begin{proof}
    First we use $w \varphi^{q}$ as a test function in \eqref{localeqn}, and obtain:
    \begin{equation}\label{eqn1}
        \begin{aligned}
            \int_{\Omega} |\nabla u|^{p-2}\nabla u \nabla (w \varphi^{q}) + A_{q}(u, w \varphi^{q}, \RR^{2N})= \int_{\Omega} f w \varphi^{q}.
        \end{aligned}
    \end{equation}
Now,
    \begin{equation*}
        \begin{aligned}
             \int_{\Omega} |\nabla u|^{p-2}\nabla u \nabla (w \varphi^{q}) =& \int_{B_{r}} |\nabla u|^{p-2} \nabla u \nabla w \varphi^{q} + \int_{B_{r}} |\nabla u|^{p-2} \nabla u w q \varphi^{q-1}\nabla \varphi
             \\\geq & \int_{B_{r}} |\nabla w|^{p} \varphi^{q} - \int_{B_{r}} |\nabla w|^{p-1}   \varphi^{q-1} w|\nabla \varphi|
             \\\geq & \int_{B_{r}} |\nabla w|^{p} \varphi^{q} - \varepsilon\int_{B_{r}} (|\nabla w|^{p-1}   \varphi^{q-1})^\frac{p}{p-1} - C(\varepsilon) \int_{B_{r}} (w |\nabla \varphi|)^{p}
             \\
             \geq & \int_{B_{r}} |\nabla w|^{p}\varphi^{q}(1-\varepsilon \varphi^{\frac{(q-1)p}{p-1}-q}) - C(\varepsilon) \int_{B_{r}} (w |\nabla \varphi|)^{p}
             \\
             =&  \int_{B_{r}} |\nabla w|^{p}\varphi^{q}(1-\varepsilon \varphi^{\frac{(q-p)}{p-1}}) - C(\varepsilon) \int_{B_{r}} (w |\nabla \varphi|)^{p}.
        \end{aligned}
    \end{equation*}
   Since $0\leq \varphi\leq 1,$ we choose $\varepsilon=1/2$  to get that 
   \begin{equation}\label{eqn2}
       \begin{aligned}
         \int_{\Omega} |\nabla u|^{p-2}\nabla u \nabla (w \varphi^{q})   \geq c \int_{B_{r}} |\nabla w|^{p}\varphi^{q} - C \int_{B_{r}} w^{p} |\nabla \varphi|^{p}.
       \end{aligned}
   \end{equation}
    Clearly we also have 
    \begin{equation}\label{eqn3}
        \begin{aligned}
            \int_{\Omega} f w \varphi^{q} \leq \int_{B_{r}} |f| w \varphi^{q}.
        \end{aligned}
    \end{equation}
    Observe that 
    \begin{equation}\label{domainbrk}
        \begin{aligned}
            A_{q}(u, w \varphi^{q}, \RR^{2N})=  A_{q}(u, w \varphi^{q}, B_{r}\times B_{r}) + 2 \int_{\RR^{N}\setminus B_{r}} \int_{B_{r}} \frac{[u(x)-u(y)]^{q-1} w(x) \varphi^{q}(x)}{|x-y|^{N+qs}}.
        \end{aligned}
    \end{equation}
    By \cite[Eqn (3.3), (3.4)]{GKS23} we know 
    \begin{equation*}
        \begin{aligned}
        &\Bigg\{[u(x)-u(y)]^{q-1}-[w(x)-w(y)]^{q-1}\Bigg\}(w(x)\varphi^{q}(x)-w(y)\varphi^{q}(y))\geq 0;\\
            &\int_{\RR^{N}\setminus B_{r}} \int_{B_{r}} \frac{[u(x)-u(y)]^{q-1} w(x) \varphi^{q}(x)}{|x-y|^{N+qs}} \geq -\left(\sup_{y\in \supp \varphi} \int_{\RR^{N}\setminus B_{r}} \frac{w(x)^{q-1}}{|x-y|^{N+sq}} \right) \int_{B_{r}} w \varphi^{q}.
        \end{aligned}
    \end{equation*}
    Hence we get 
    \begin{equation*}
        \begin{aligned}
           A_{q}(u, w \varphi^{q}, \RR^{2N}) \geq A_{q}(w, w \varphi^{q}, B_{r}\times B_{r}) -\left(\sup_{y\in \supp \varphi} \int_{\RR^{N}\setminus B_{r}} \frac{w(x)^{q-1}}{|x-y|^{N+sq}} \right) \int_{B_{r}} w \varphi^{q}.
        \end{aligned}
    \end{equation*}
    By  \cite[pp.10]{GKS23} we have
    \begin{equation*}
        \begin{aligned}
          \relax  [w(x)-w(y)]^{q-1} (w(x)\varphi^{q}(x)-w(y)\varphi^{q}(y)) &\geq |w(x)-w(y)|^{q}\frac{\varphi(x)^{q}+\varphi(y)^{q}}{4}
            \\
            &-C(w(x)+w(y))^{q}|\varphi(x)-\varphi(y)|^{q}.
        \end{aligned}
    \end{equation*}
    Hence we infer that
    \begin{equation}\label{eqn4}
        \begin{aligned}
     A_{q}(u, w \varphi^{q}, \RR^{2N}) \gtrsim & \int_{B_{r}} \int_{B_{r}} \frac{|w(x)-w(y)|^{q} ( \varphi^{q}(x)+ \varphi^{q}(y))}{|x-y|^{N+sq}} \\&- \int_{B_{r}} \int_{B_{r}} \frac{|\varphi(x)-\varphi(y)|^{q}(w(x)+w(y))^{q}}{|x-y|^{N+sq}}
     \\
     & -\left(\sup_{y\in \supp \varphi} \int_{\RR^{N}\setminus B_{r}} \frac{w(x)^{q-1}}{|x-y|^{N+sq}} \right) \int_{B_{r}} w \varphi^{q}.
        \end{aligned}
    \end{equation}
    Combining \eqref{eqn1}-\eqref{eqn4},  we get our desired result.
\end{proof}
\begin{remark}
   For \( q \leq p \), a similar Caccioppoli inequality can be established by employing $w\varphi^{p}$ as test function. Furthermore, following the approach outlined in \cite[Proposition 3.2]{GKS23}, one can derive the local boundedness result. However, we do not provide a detailed proof here, as we rely on the interior regularity framework presented in \cite{FM22}.
\end{remark}
In a manner similar to \cite[Lemma 2.4]{BOS22}, we obtain the following lemma:
\begin{lemma}\label{alglem}
Let $1<p\leq q<\infty $ and $w\in \wsqq(B_{r}).$ Then 
\begin{equation*}
    \begin{aligned}
        \intavg_{B_{r}} |\frac{w}{r^{s}}|^{q} \leq C \left[\left(\frac{|\supp w|}{|B_{r}|}\right)^{\frac{sq}{N}}  \intavg_{B_{r}} \int_{B_{r}} \frac{|w(x)-w(y)|^{q}}{|x-y|^{N+sq}}+ \left(\frac{|\supp w|}{|B_{r}|}\right)^{q-1}  \intavg_{B_{r}}\left|\frac{w}{r^{s}}\right|^{q}\right]
    \end{aligned}
\end{equation*}
for some $c(N,s,q)>0.$
\end{lemma}
We now prove an important step towards proving the local H\"older regularity result.
\begin{theorem}\label{localbdd}
    Let $1\leq p \leq q<\infty,$ and $u\in\mathcal{X}(\Omega) $ be a local solution to \eqref{localeqn}. If $f\in L^{\infty}_{\loc}(\Omega)$ then $u \in L^{\infty}_{\loc}(\Omega).$
\end{theorem}
\begin{proof}
    Let $B_{r}\equiv B_{r}(x_{0})\subset B_{2r}(x_0)\subset\subset\Omega$ be a fixed ball with $r\leq 1.$ For $\frac{r}{2}\leq \varrho<\sigma\leq r$ and $k>0,$ we denote $A^{+}(k,\varrho)=\{x\in B_{\varrho}: u(x)\geq k \}.$ For $w=(u-k)_{+}$ in Lemma \ref{alglem} we obtain
    \begin{equation}\label{ccp1}
    \begin{aligned}
     \varrho^{-sq}   \intavg_{B_{\varrho}} {(u-k)_{+}}^{q} \leq & C\left[\left(\frac{|A^{+}(k,\varrho)|}{|B_{\varrho}|}\right)^{\frac{sq}{N}}  \intavg_{B_{\varrho}} \int_{B_{\varrho}} \frac{|(u-k)_{+}(x)-(u-k)_{+}(y)|^{q}}{|x-y|^{N+sq}}\right.
     \\
     &\hspace{9em}+\left.\left(\frac{|A^{+}(k,\varrho)|}{|B_{\varrho}|}\right)^{q-1}  \intavg_{B_{\varrho}} (u-k)_{+}^{q}\varrho^{-sq}\right]
    \end{aligned}
\end{equation}
for $C(N,s,q)>0.$
We now fix $0<h<k$ and note that for $x\in A^{+}(k,\varrho) \subset A^{+}(h,\varrho)$ we get $(u-h)_{+}\geq k-h$ and $(u-h)_{+}\geq (u-k)_{+}.$ This gives 
\begin{equation}\label{ccp2}
    \begin{aligned}
      \frac{  |A^{+}(k,\varrho)|}{|B_{\varrho}|} &\leq \frac{1}{|B_{\varrho}|}\int_{A^{+}(k,\varrho)} \frac{(u-h)_{+}^{q}}{(k-h)^{q}} \leq \frac{1}{(k-h)^{q}} \left(\frac{\sigma}{\varrho}\right)^{N} \intavg_{B_{\sigma}} (u-h)_{+}^{q}\text{ and }
      \\
      &\int_{B_{\sigma}} (u-k)_{+} \leq \frac{1}{(k-h)^{q-1}} \int_{B_{\sigma}}  (u-h)_{+}^{q} .
    \end{aligned}
\end{equation}
We then choose a cut off function $\varphi\in C_{0}^{\infty}(B_{\frac{\varrho+\sigma}{2}})$ satisfying $0\leq \varphi\leq 1,$ $\varphi\equiv 1$ in $B_{\varrho}$ and $|D \varphi|\leq \frac{4}{\sigma-\varrho}.$ Lemma \ref{cacciopoli} gives us 
\begin{equation}\label{ccp3}
    \begin{aligned}
      \intavg_{B_{\varrho}} \int_{B_{\varrho}}\frac{|w(x)-w(y)|^{q}}{|x-y|^{N+sq}} &\leq \frac{1}{|B_{\varrho}|}\left(\int_{B_{\varrho}}|Dw|^{p}+ \int_{B_{\varrho}} \int_{B_{\varrho}}\frac{|w(x)-w(y)|^{q}}{|x-y|^{N+sq}}\right) 
            \\
            &\leq \frac{C}{|B_{\varrho}|}\left[ \int_{B_{\sigma}} w^{p} |D\varphi|^{p}+\int_{B_{\sigma}} \int_{B_{\sigma}} \frac{|\varphi(x)-\varphi(y)|^{q}(w(x)+w(y))^{q}}{|x-y|^{N+sq}} \right.
            \\& \hspace{3.5em}+ \left.\left(\sup_{x\in \supp \varphi} \int_{\RR^{N}\setminus B_{\sigma}} \frac{w(y)^{q-1}\thinspace dy}{|x-y|^{N+sq}}\right) \int_{B_{\sigma}} w \varphi^{q} +\int_{B_{\sigma}}|f|w\varphi^{q}\right]
    \end{aligned}
\end{equation}
for $C(N,s,p,q)>0.$ Next observe that for $x\in \supp \varphi$ and $y \in \RR^{N}\setminus B_{\sigma},$
\begin{equation*}
    \begin{aligned}
        \frac{|y-x_0|}{|y-x|}\leq 1 + \frac{|x-x_0|}{|y-x|}\leq 1+ \frac{\sigma+\varrho}{\sigma-\varrho} \leq 2 \frac{\sigma+\varrho}{\sigma-\varrho}.
    \end{aligned}
\end{equation*}
Then we derive:
\begin{equation*}
    \begin{aligned}
        \intavg_{B_{\varrho}} \int_{B_{\varrho}}\frac{|w(x)-w(y)|^{q}}{|x-y|^{N+sq}} &\lesssim \left(\frac{\sigma}{\varrho}\right)^{N}\Bigg[\frac{1}{(\sigma-\varrho)^{p}}\intavg_{B_{\sigma}} w^{p}+ \frac{1}{(\sigma-\varrho)^{q}}\intavg_{B_{\sigma}} (u-h)_{+}^{q}\int_{B_{\sigma}} \frac{dx\thinspace dy}{|x-y|^{N+(s-1)q}}
         \\
         &
         \hspace{3.0em}+\left(\sup_{x\in \supp \varphi} \int_{\RR^{N}\setminus B_{\sigma}} \frac{w(y)^{q-1}\thinspace dy}{|x-y|^{N+sq}} \right) \intavg_{B_{\sigma}} w  +\|f\|_{\infty}\intavg_{B_{\sigma}}w \Bigg],
         \\
         &\lesssim \left(\frac{\sigma}{\varrho}\right)^{N}\Bigg[\frac{1}{(\sigma-\varrho)^{p}}\intavg_{B_{\sigma}} w^{p}+ \frac{\varrho^{(1-s)q}}{(\sigma-\varrho)^{q}}\intavg_{B_{\sigma}} (u-h)_{+}^{q}+\|f\|_{\infty}\intavg_{B_{\sigma}}w
         \\
         &
        \hspace{3.0em} +\left(\frac{\sigma+\varrho}{\sigma-\varrho}\right)^{N+sq}\left(\int_{\RR^{N}\setminus B_{\sigma}} \frac{w(y)^{q-1}\thinspace dy}{|y-\xo|^{N+sq}} \right) \intavg_{B_{\sigma}} w \Bigg]
    \end{aligned}
\end{equation*}
in a similar manner as presented in \cite[pp. 126]{BOS22}. Hence
\begin{equation}\label{ccp4}
    \begin{aligned}
     \intavg_{B_{\varrho}} \int_{B_{\varrho}}\frac{|w(x)-w(y)|^{q}}{|x-y|^{N+sq}}     &\lesssim \left(\frac{\sigma}{\varrho}\right)^{N}\Bigg[\frac{1}{(\sigma-\varrho)^{p}}\intavg_{B_{\sigma}} w^{p}+ \frac{r^{(1-s)q}}{(\sigma-\varrho)^{q}}\intavg_{B_{\sigma}} (u-h)_{+}^{q}+\|f\|_{\infty}\intavg_{B_{\sigma}}w
         \\
         &
       \hspace{4.0em}  +\left(\frac{r}{\sigma-\varrho}\right)^{N+sq}\left(\int_{\RR^{N}\setminus B_{\sigma}} \frac{w(y)^{q-1}\thinspace dy}{|y-\xo|^{N+sq}} \right) \intavg_{B_{\sigma}} w \Bigg].
    \end{aligned}
\end{equation}
Combining \eqref{ccp1},\eqref{ccp4} yields:
\begin{equation*}
    \begin{aligned}
        \varrho^{-sq}\intavg_{B_{\varrho}} (u-k)_{+}^{q} &\lesssim \left(\frac{\sigma}{\varrho}\right)^{N} \left(\frac{|A^{+}(k,\varrho)|}{|B_{\varrho}|}\right)^{\frac{sq}{N}} \Bigg[\frac{1}{(\sigma-\varrho)^{p}}\intavg_{B_{\sigma}} w^{p}+ \frac{r^{(1-s)q}}{(\sigma-\varrho)^{q}}\intavg_{B_{\sigma}} (u-h)_{+}^{q}
        \\&\hspace{4.5cm}+\|f\|_{\infty}\intavg_{B_{\sigma}}w
         \\
         &
        \hspace{12em} +(\frac{r}{\sigma-\varrho})^{N+sq}\left(\int_{\RR^{N}\setminus B_{\sigma}} \frac{w(y)^{q-1}\thinspace dy}{|y-\xo|^{N+sq}} \right) \intavg_{B_{\sigma}} w \Bigg]
         \\
         & \hspace{12em}+ \left(\frac{|A^{+}(k,\varrho)|}{|B_{\varrho}|}\right)^{q-1}  \intavg_{B_{\varrho}} (u-k)_{+}^{q}\varrho^{-sq}.
    \end{aligned}
\end{equation*}
Since $p\leq q,$ we observe that 
\begin{equation*}
    \begin{aligned}
          (u-k)_{+}^{p}&\leq  \frac{(u-h)_{+}^{q}}{(k-h)^{q-p}} ,
        \\
        \int_{B_{\sigma}} (u-k)_{+}&\leq \frac{1}{(k-h)^{q-1}} \int_{B_{\sigma}} (u-h)^{q}_{+}.
    \end{aligned}
\end{equation*}
Combining this with \eqref{ccp2}, we obtain
\begin{equation*}
    \begin{aligned}
        &\varrho^{-sq}\intavg_{B_{\varrho}} (u-k)_{+}^{q} \lesssim \left(\frac{\sigma}{\varrho}\right)^{N(\frac{sq}{N}+1)}\frac{1}{(k-h)^{\frac{sq^{2}}{N}+(q-p)}(\sigma-\varrho)^{p}} \left(\intavg_{B_{\sigma}}(u-h)_{+}^{q}\right)^{1+\frac{sq}{N}}
         \\
         &
         + \left(\frac{\sigma}{\varrho}\right)^{N(\frac{sq}{N}+1)}\frac{r^{(1-s)q}}{(k-h)^{\frac{sq^{2}}{N}}(\sigma-\varrho)^{q}} \left(\intavg_{B_{\sigma}} (u-h)_{+}^{q}\right)^{1+\frac{sq}{N}}
         \\
         & + \left(\frac{\sigma}{\varrho}\right)^{N(\frac{sq}{N}+1)} \frac{\|f\|_{L^{\infty}} }{(k-h)^{\frac{sq^{2}}{N}+q-1}} \left(\intavg_{B_{\sigma}}(u-h)_{+}^{q}\right)^{1+\frac{sq}{N}}
         \\
         & + \left(\frac{\sigma}{\varrho}\right)^{N(\frac{sq}{N}+1)} \frac{r^{N+sq}}{(\sigma-\varrho)^{N+sq}(k-h)^{\frac{sq^2}{N}+q-1}}\left(\int_{\RR^{N}\setminus B_{\sigma}} \frac{(u-k)_{+}^{q-1}\thinspace dy}{|y-\xo|^{N+sq}}\right)\left(\intavg_{B_{\sigma}}(u-h)_{+}^{q}\right)^{1+\frac{sq}{N}}
         \\
         &
         + \left(\frac{\sigma}{\varrho}\right)^{N(q-1)}\frac{r^{-sq}}{(k-h)^{q(q-1)}}\left(\intavg_{B_{\sigma}}(u-h)_{+}^{q}\right)^{q}.
    \end{aligned}
\end{equation*}
Letting $\alpha = \max\{N(q-1), N + sq\}$ and using $\sigma/\varrho > 1$, we find that
\begin{equation*}
    \begin{aligned}
   &\varrho^{-sq}\intavg_{B_{\varrho}} (u-k)_{+}^{q} \lesssim \left(\frac{\sigma}{\varrho}\right)^{\alpha}\Bigg[\frac{1}{(k-h)^{\frac{sq^{2}}{N}+{q-p}}(\sigma-\varrho)^{p}} \left(\intavg_{B_{\sigma}}(u-h)_{+}^{q}\right)^{1+\frac{sq}{N}}
   \\&
   +\frac{r^{(1-s)q}}{(k-h)^{\frac{sq^{2}}{N}}(\sigma-\varrho)^{q}} \left(\intavg_{B_{\sigma}} (u-h)_{+}^{q}\right)^{1+\frac{sq}{N}}
   +\frac{\|f\|_{L^{\infty}} }{(k-h)^{\frac{sq^{2}}{N}+q-1}} \left(\intavg_{B_{\sigma}}(u-h)_{+}^{q}\right)^{1+\frac{sq}{N}}
   \\&
   +\frac{r^{N+sq}}{(\sigma-\varrho)^{N+sq}(k-h)^{\frac{sq^2}{N}+q-1}}\left(\int_{\RR^{N}\setminus B_{\sigma}} \frac{(u-k)_{+}^{q-1}}{|x-\xo|^{N+sq}}\right)\left(\intavg_{B_{\sigma}}(u-h)_{+}^{q}\right)^{1+\frac{sq}{N}}
  \\&
  + \frac{r^{-sq}}{(k-h)^{q(q-1)}}\left(\intavg_{B_{\sigma}}(u-h)_{+}^{q}\right)^{q}
   \Bigg].    
    \end{aligned}
\end{equation*}
For $i=0,1,2 \ldots,$ set 

\begin{minipage}{.5\linewidth}
\begin{equation*}
  \begin{aligned}
 k_{0}&>1,
      \\
      k_{i}&:= 2k_{0}(1-2^{-i-1}),
  \end{aligned}
\end{equation*}
\end{minipage}%
\begin{minipage}{.5\linewidth}
\begin{equation*}
  \begin{aligned}
     \sigma_{i}&=\frac{r}{2}(1+2^{-i}),
      \\
     y_{i}&:=\int_{A_{+}(k_{i},\sigma_{i})}(u-k_{i})_{+}^{q}.
  \end{aligned}
\end{equation*}
\end{minipage}
Notably since $u\in L^{q}_{sq}(\RR^N) $ we have 
\begin{equation*}
    \begin{aligned}
        \int_{\RR^{N}\setminus B_{\sigma_{i}}}\frac{(u-k)_{+}^{q-1}}{|x-\xo|^{N+sq}}\leq\int_{\RR^{N}\setminus B_{r/2}}\frac{(u-k)_{+}^{q}}{|x-\xo|^{N+sq}}<\infty.
    \end{aligned}
\end{equation*}
\\
Moreover, since $u^{q}\in L^{1}_{\loc}(\Omega)$ we conclude
\begin{equation*}
    \begin{aligned}
        y_{0}=\int_{A^{+}(k_{0},r)}(u-k_{0})_{+}^{q}\rightarrow 0 \text{ as } k_{0}\rightarrow \infty.
    \end{aligned}
\end{equation*}
We first consider $k_0>1$ sufficiently large such that $y_{i}\leq y_{i-1}\leq\ldots\leq y_{0}\leq 1.$ Then
\begin{equation*}
    \begin{aligned}
        y_{i+1}\lesssim [2^{i(q+\frac{sq^2}{N})}y_{i}^{1+\frac{sq}{N}}+ 2^{i(q+\frac{sq^{2}}{N})}y_{i}^{1+\frac{sq}{N}}
        +2^{i(N+sq+\frac{sq^{2}}{N}+q-1)}y_{i}^{1+\frac{sq}{N}}
        +2^{i(\frac{sq^{2}}{N}+q-1)}y_{i}^{1+\frac{sq}{N}}
        +2^{iq(q-1)}y_{i}^{q}].
    \end{aligned}
\end{equation*}
Choose $\beta=\min\{\frac{sq}{N},q-1\}$ and $\theta=q^{2}+N$ to get that 
\begin{equation*}
    \begin{aligned}
        y_{i+1}\leq C 2^{i\theta}y_{i}^{1+\beta}
    \end{aligned}
\end{equation*}
for some $C(N,s,p,q)>0.$ Finally we can choose $k_0>1$ so large that $y_{0}\leq C^{\frac{-1}{\beta}}2^{\frac{-\theta}{\beta^{2}}}.$ Then using Lemma \ref{itrlrm} we get that $y_{i}\rightarrow0$ as $i\rightarrow \infty$ which means that $u\leq 2k_{0}$ in $B_{r/2}.$ Applying the same argument to $-u$ we consequently obtain $u\in L^{\infty}(B_{r/2}).$
\end{proof}

\begin{remark}
    From the proof, it is evident that the local boundedness result holds if we assume  $u \in \wpp_{\loc}(\Omega)\cap\wsqq_{\loc}(\Omega)\cap L^{q-1}_{sq}(\RR^N)$ instead of $u \in \mathcal{X}(\Omega).$
\end{remark}

\subsection{Interior Regularity}

First, we will rewrite \cite[Theorem 6]{FM22} in the light of \cite[Remark 3]{FM22} and zero Dirichlet boundary condition for the operator $\fp + \fqs.$
\begin{theorem}(Theorem 6, \cite{FM22})\label{intfm22}
    Let $u\in \wpp_{\loc}(\Omega)\cap\wsqq_{\loc}(\Omega)$ solves \eqref{localeqn}and 
    \begin{equation*}
        \begin{aligned}
            \int_{\RR^{N}}\left[\frac{|u(x)|^{q}}{(1+|x|)^{N+sq}}\thinspace dx\right]^{\frac{1}{q}}<\infty.
        \end{aligned}
    \end{equation*}
    If $u \in L^{\infty}_{\loc}(\Omega),$  $f\in L^{d}_{\loc}(\Omega)$ for some $d>N$ and $q>p,$ then $u \in C^{1,\alpha}_{\loc}(\Omega)$ for some $\alpha\in(0,1).$
\end{theorem}
We will now give the proof of our first main result : local H\"older regularity. \\
\textit{\underline{Proof of Theorem \ref{interiorregularity} }}
\\[3mm]
We will prove this result under two conditions: 
$ p<q$ or $ q\leq p.$\\
\underline{Case (a):} $1<p<q<\infty.$  Given $f\in L^\infty_{loc}(\Omega),$ from Lemma \ref{localbdd} we have $u\in L^\infty_{loc}(\Omega).$ Now, Theorem \ref{intfm22} is applicable and $u\in C^{\alpha}_{loc}(\Omega)$ for each $\alpha\in (0,1).$ Indeed, $u\in C^{1,\alpha'}_{loc}(0,1)$ for some $\alpha'\in (0,1).$\\[3mm]
\underline{Case (b):} $1<q\leq p <\infty.$  This is an easier case where it is not necessary to prove an $L^\infty$ regularity before obtaining the H\"older type results. This is because, H\"older regularity result in this case is direct from  Theorem 3, Theorem 5, Remark 3 and Section 7 of Mingione and D. Filippis\cite{FM22}. When $q\leq p,$ by definition, we consider a local solution $u\in \mathcal{X}(\Omega).$ From Remark 3 of \cite{FM22}, if $f\in L^d(\Omega),$ one can prove that $u$ belongs to $C^{\alpha}_{loc}(\Omega)$ directly without proving $L^\infty_{loc}$ regularity. Moreover, the same remark  mentions that the $C^{\alpha}(\Omega')$ norm of $u$ depends only on the data given by  
\begin{equation*}
    \begin{aligned}
        &\text{data}_{h} := \left(n, p, s, q, \Omega, \|f\|_{L^N(\Omega')},\|u\|_{L^p(\Omega')}, \|u\|_{L^{q}(\RR^N)}\right) 
    \end{aligned}
\end{equation*}
while $\|u\|_{L^q(\RR^N)}$ can be replaced with $\|u\|_{L^q_\loc(\Omega')}$ plus the tail space norm of $u$ in $L^{q}_{sq}(\RR^N).$ Since $u\in L^{q}_{sq}(\RR^N),$ thus the snail term defined in (3.2) of \cite{FM22} is well defined. Now invoking results mentioned in section 7 of \cite{FM22}, we can consider $f\in L^d_{loc}(\Omega)$ and prove the results analogous to Theorem 3 and Theorem 5 of \cite{FM22}, i.e $u\in C^{\alpha}_{loc}(\Omega)$ for all $\alpha\in (0,1)$ and $u\in C^{1,\theta}_{loc}(\Omega)$ for some $\theta\in (0,1).$ 
This completes the proof of Theorem \ref{interiorregularity}. \hfill\qed 

\subsection{Boundary Regularity}

In this section we consider the equation 
\begin{equation}\label{singregeqn}
    \begin{aligned}
        \fp u +\fqs u&= f(x)\text{ in } \Omega,
        \\
        u&=0 \text{ in } \Omega^{c}
    \end{aligned}
\end{equation}
where $f\in L^{\infty}_{\loc}(\Omega)$ such that   \begin{equation*}
0<\Theta \leq f(x) \leq C_{f}  \dis^{-\gamma}(x) \text{ in } \Omega
\end{equation*}
for some $\Theta, C_{f}>0$ and we fix the constants in this section. From the previous results, we have already established that the solution $u$ is locally Hölder continuous. Next, our goal is to analyze the boundary behavior of the solution under the  case $0<\gamma<1.$ 
\subsubsection{Boundary H\"older Regularity }
\begin{theorem}\label{hldrtm}
Let $0<\gamma<1$ and $u\in \wop(\Omega)\cap\wsq(\Omega) $ be a weak solution of 
 \eqref{singregeqn}. Then,
 \begin{enumerate}
     \item[(a)]  for every $\mu\in (\frac{qs}{q-1},1,)$  $C_1 \dis(x) \leq u(x) \leq    C_2 \dis^{\mu}(x)$ for all $x\in \Omega.$
     \item[(b)]  Additionally, $u\in C^{\mu}(\overline{\Omega})$ for each $\mu\in (0,1).$
 \end{enumerate} 
\end{theorem}
\begin{proof}
    Consider the unique solution $w_{\Theta}\in \wop(\Omega)\cap\wsq(\Omega)\cap C^{1,\alpha}(\Bar{\Omega})$ of the following equation 
    \begin{equation*}
        \begin{aligned}
            \fp w_{\Theta}+ \fqs w_{\Theta} &= \Theta,\quad w_{\Theta}>0 \text{ in } \Omega,
            \\
            w_{\Theta}&=0 \text{ on } \Omega^{c}.
        \end{aligned}
    \end{equation*}
By the weak comparison principles we get  $w_\Theta\leq u$ in $\Omega.$ Consequently by the Hopf's Lemma \cite[Theorem 1.2]{AC23}  we obtain 
\begin{equation*}
    C_1 \dis \leq u \text{ in } \Omega.
\end{equation*}
For $\rho>0$  we define the extension of distance function  $\dis_e(x)$ in a way consistent with the nonlocal case \cite[pp. 23]{GKS23} as follows:
\begin{equation}\label{dedef}
    \dis_{e}(x)= \left\{
    \begin{aligned}
        \dis(x) &\text{ if } x \in \Omega,\\
        -\dis(x) &\text{ if } x \in (\Omega^{c})_{\rho},\\
        -\rho &\text{ otherwise }
    \end{aligned}
    \right. 
\end{equation}
where $(\Omega^{c})_{\rho}:=\{x\in \Omega^c: \dis(x,\partial\Omega)<\rho\}.$
Additionally, for $\alpha\in (\frac{qs}{q-1},1),$  we set
 \begin{equation*}
\Bar{\omega}_\rho(x)=\left\{
\begin{aligned}
    &(\dis_{e}(x))_{+}^{\alpha} \text{ in } \Omega \cup (\Omega^{c})_{\rho},
    \\
    &0 \text{ otherwise.}
\end{aligned}\right.   
\end{equation*}
Using \cite[Eqn (5.23 )]{GKS24} since $\alpha>\frac{qs}{q-1},$ we find that for sufficiently small $\rho>0$
\begin{equation*}
    \begin{aligned}
         \fqs \Bar{\omega}_\rho= h \text{ in } \Omega_\rho \text{ where } h \in L^{\infty}(\Omega_\rho). 
    \end{aligned}
\end{equation*}
Furthermore, for $\rho$ small enough, we have that $\Delta \dis\in L^{\infty}(\Omega_\rho)$ and $|\nabla\dis|=1.$ Thus for any $\psi\in C_{c}^{\infty}(\Omega_\rho)$ with $\psi\geq 0$ and  we can deduce that 
\begin{equation*}
    \begin{aligned}
        \int_{\Omega_{\rho}} \fp \Bar{\omega}_\rho \psi &= \alpha^{p-1}   \int_{\Omega_{\rho}} [-\Delta d(\dis(x))^{(\alpha-1)(p-1) }\psi + (p-1)(1-\alpha)(\dis(x))^{(\alpha p -\alpha-p) }\psi].
    \end{aligned}
\end{equation*}
Therefore, 
\begin{equation}\label{bdbehaviourbarrier}
    \begin{aligned}
    \fp \Bar{\omega}_\rho      \geq C (\dis(x))^{\alpha p - \alpha -p} \text{ in } \Omega_\rho \text{ for } \rho \text{ small enough.}
    \end{aligned}
\end{equation}
Since $\alpha<1<\frac{p-\gamma}{p-1},$ we deduce that
\begin{equation*}
    \begin{aligned}
        \alpha p -\alpha-p<-\gamma.
    \end{aligned}
\end{equation*}
Hence for $\rho$ sufficiently small we deduce that 
\begin{equation}
   \begin{aligned}
        \fp (\Gamma \wro) + \fqs (\Gamma \wro) &\geq C \Gamma^{p-1} (\dis(x))^{\alpha p - \alpha -p} -\Gamma^{q-1} \|h\|_{L^{\infty}(\Omega_\rho)}\\
    &\geq \Gamma^{p-1} \frac{C}{2} (\dis(x))^{\alpha p - \alpha -p}\\
    &\geq  C_{f}\dis(x)^{-\gamma}  \text{ in } \Omega_{\rho}
   \end{aligned}
\end{equation}
provided that $ \Gamma^{p-1}\frac{C}{2}\geq C_{f}.$
So we can choose $\Gamma$ large enough such that $\wro$ is supersolution to \eqref{singregeqn} in the weak sense in $\Omega_\rho.$ Now observe that since $f(x)\in L^{\infty}_{\loc}(\Omega)$ using \eqref{interiorholderreg} of Theorem \ref{interiorregularity} we have $u\in C^{\alpha}_{\loc}(\Omega)$ for each $\alpha\in(0,1).$ Now by choosing $\Gamma$ larger if required 
\begin{equation*}
    \begin{aligned}
        u\leq \|u\|_{L^{\infty}(\Omega\setminus\Omega_{\rho})}\leq C_{\rho} \leq \Gamma\dis^{\alpha}=\Gamma\wro \text{ in } \Omega\setminus\Omega_{\rho}.
    \end{aligned}
\end{equation*}
Using the weak comparison principle in $\Omega_\rho$ we prove part $(a)$ of the theorem. Part (b) follows from Part $(a)$ and the interior H\"older regularity result since $\Omega$ has smooth boundary. 
\end{proof} 

\subsubsection{Boundary Gradient H\"older Regularity}
In this subsection, we aim to prove the boundary H\"older regularity for \eqref{singregeqn}. We adapt the approach presented in \cite[Theorem 1.1]{AC23} but with the singular data and also need to revise the method outlined in \cite[Theorem 4]{FM22}. We begin by flattening the boundary. Following the discussion in \cite[pp. 7-8]{AC23}, we locally straighten the boundary around any point  $x_0 \in \partial\Omega.$ Then, based on the argument in \cite[Section 5]{FM22}, we establish the existence of a global  $C^{1,\alpha}$-diffeomorphism  $\mathcal{T}$  on  $\RR^N$, such that
\begin{equation*}
    \begin{aligned}
        &\T(\xo)=\xo,\\
    & B^{+}_{r_{0}}(\xo) \subset \T(\Omega_{3r_{0}}(\xo)) \subset B^{+}_{4r_{0}}(\xo),
    \\& \Gamma_{r_{0}}(\xo) \subset \T(\partial\Omega\cap B_{3r_{0}}(\xo)) \subset \Gamma_{4r_{0}}(\xo) \text{ for } r_0\in (0,1]
    \end{aligned}
\end{equation*}
where 
\begin{equation*}
    \begin{aligned}
        \Omega_{r}(\xo):=\Omega\cap B_{r}(\xo)\quad\text{ and }\quad \Gamma_{r}(\xo):= B_{r}(\xo)\cap \{x_n=0\}.
    \end{aligned}
\end{equation*}
Write $\Sl:=\T^{-1}$ and $\mathfrak{c}:=|\mathcal{J_S}|,$ where $\mathcal{J_S}$ denotes the Jacobian determinant of the inverse $\Sl.$ Let $\widetilde{\Omega}:=\T(\Omega), \wkg:=\ck(\kg\circ\Sl),\wf:=f\circ{\T}^{-1}$ and $\widetilde{u}:=u\circ\Sl.$ 
For $y\in\widetilde{\Omega }\cap B^{+}_{r_{0}}(\xo)$ , we get that $C_2 y_N \leq \dis(S(y))\leq C_1 y_N.$ Then it is easy to see that $\wf\in L^{\infty}_{\loc}(\widetilde{\Omega})$ such that 
\begin{equation*}
    \begin{aligned}
        \frac{C_1}{(y_{N})^{\gamma}} \leq \wf(y) \leq \frac{C_2}{(y_{N})^{\gamma}}.
    \end{aligned}
\end{equation*}
We have $\wu$ weakly solves 
\begin{equation}\label{wueqn}
    \begin{aligned}
        -\Div \widetilde{A}(\cdot,D\wu)+\widetilde{Q}_N\wu&= \wf,\quad \wu>0 \text{ in } \wO,\\
        \wu&=0 \text{ in } \wO^{c}
    \end{aligned}
\end{equation}
where 
\begin{equation*}
    \begin{aligned}
        \widetilde{A}(x,z)&:=\mathfrak{c}(x) A(\Sl(x),z(D\T\circ \Sl)(x))(D\T\circ \Sl)(x)^{T},
        \\
        A(x,z)&:=|z|^{p-2}z,
        \\
        \widetilde{Q}_N u(x)&:= 2 \pv \int_{\RR^N} |u(x)-u(y)|^{q-1} \widetilde{K}(x,y) \thinspace dy \text{ with }
        \\
        \widetilde{K}(x,y)&:= \mathfrak{c}(x)\mathfrak{c}(y)\frac{1}{|\Sl(x)-\Sl(y)|^{N+sq}}.
    \end{aligned}
\end{equation*}
Consider any point $\widetilde{x}_0 \in \Gamma_{r_{0/2}}(x_0)$ and radius  $\rho \in (0,r_{0}/4]$. Since the data $\widetilde{f}$ is integrable in \cite{FM22}, the authors can utilize the H\"older's inequality to estimate the term involving this data. At this stage, we are not addressing the integral involving the locally integrable singular data near the boundary. However, by carefully adapting the proof of \cite[Lemma 5.1]{FM22}, we can derive the following result.

\begin{lemma} \label{ccplem}
   Let $\wu$ be a solution of \eqref{wueqn} in the sense of Definition \ref{soldef}. Then 
   \begin{equation*}
       \begin{aligned}
           \dashint_{B^{+}_{\rho/2}(\widetilde{x}_0)} |D\tu|^{p}\thinspace dx & +  \int_{B_{\rho/2}(\widetilde{x}_0)}  \dashint_{B^{+}_{\rho/2}(\widetilde{x}_0)} \frac{|\tu(x)-\tu(y)|^q}{|x-y|^{N+sq}} \thinspace dx \thinspace dy
        \\& \leq c\rho^{-p} \dashint_{B^{+}_{\rho}(\widetilde{x}_0)} |\tu|^{p} + c \int_{\RR^N\setminus B_{\rho}(\widetilde{x}_0)}\frac{|\tu(x)-\tu(y)|^{q}}{|x-y|^{N+sq}}
        + \dashint_{B^{+}_{\rho} (\widetilde{x}_0)} \eta^{m} \ck(x)\wf\wu 
       \end{aligned}
   \end{equation*}
   where $m=\max\{p,q\};$ and $\eta\in C_{0}^{1}(B_{\tau_{2}})$ such that $\rho/2\leq \tau_{1}\leq \tau_{2}\leq \rho,$ $1_{B_{\tau_{1}}}\leq \eta \leq 1_{B_{(3\tau_{2}+\tau_{1})/4}}$ with $|D\eta|\leq C/(\tau_{2}-\tau_{1})$ as defined in \cite[Lemma 5.1]{FM22}.
\end{lemma}
We now present the following preliminary estimates. Throughout this subsection, unless stated otherwise, the constants will depend only on  $N, p, q, s, \Omega, \|f\|_{L^\infty}$, and  $\beta$, where  $u \in C^{\beta}(\Bar{\Omega})$.
\begin{lemma}\label{ulem}
   Let $\wu$ be a solution of \eqref{wueqn}. Then $\wu\in C^{\beta}(\RR^N)$ for each $\beta\in(0,1)$ i.e. $\|\wu\|_{C^{\beta}(\RR^N)} \leq C_\beta.$ Moreover, it satisfies 
   \begin{equation}\label{uleme1}
       \begin{aligned}
           \int_{B_{t}(\widetilde{x}_0)}  \dashint_{B^{+}_{t}(\widetilde{x}_0)} \frac{|\tu(x)-\tu(y)|^q}{|x-y|^{N+sq}} \thinspace dx \thinspace dy \leq C_\beta t^{(\beta-s)q} \text{ for all } s<\beta<1 \text{ and } t\in (0,r_0/4);
       \end{aligned}
   \end{equation}
      \begin{equation}\label{uleme2}
       \begin{aligned}
        \int_{\RR^N\setminus B_{t}(\widetilde{x}_0) } \frac{|\tu(y)-(\tu)_{B_{\rho}(\widetilde{x}_0)}|^q}{|y-x_0|^{N+sq}}\thinspace dy \leq C \text{ for all }  t\in (0,r_0/4);
       \end{aligned}
   \end{equation}  
   \begin{equation}\label{uleme3}
       \begin{aligned}
           \dashint_{B^{+}_{\rho/2}(\widetilde{x}_0)} |D\tu|^{p}\thinspace dx \leq C_\lambda \rho^{-\lambda p} +\dashint_{B^{+}_{\rho} (\widetilde{x}_0)} \eta^{m} \ck(x)\wf\wu ,
           \text{ for all } \lambda<1-s
\text{ and }  \rho\in (0,r_0/4).
       \end{aligned}
   \end{equation}
\end{lemma}
\begin{proof}
From Theorem \ref{hldrtm}, we deduce that $\tu \in C^{\beta}(\RR^N)$ for each $0<\beta<1.$ This regularity will now be utilized to establish \eqref{uleme1}. Consider the following estimate:
\begin{equation*}
    \begin{aligned}
      \int_{B_{t}(\widetilde{x}_0)}  \dashint_{B^{+}_{t}(\widetilde{x}_0)} \frac{|\tu(x)-\tu(y)|^q}{|x-y|^{N+sq}} \thinspace dx \thinspace dy &\leq \relax[\wu]^{q}_{C^{\beta}(\RR^N)}    \int_{B_{t}(\widetilde{x}_0)}  \dashint_{B^{+}_{t}(\widetilde{x}_0)}  \frac{dx\thinspace dy}{|x-y|^{N+sq}}
      \\
      &\leq C_\beta \int_{B_{2t}(\widetilde{x}_0)} \frac{dz}{|z|^{n+(s-\beta)q}} \leq C_\beta t^{(\beta-s)q}
    \end{aligned}
\end{equation*}
where we used the change of variable $z=x-y$ and used the inclusion  $B_{t}(\widetilde{x}_0)-y\subset B_{2t}(\widetilde{x}_0) $ for each $y\in B_{t}(\widetilde{x}_0).$
\\
A similar calculation to \cite[pp. 9]{AC23} shows that 
\begin{equation}
    \begin{aligned}
        \int_{\RR^N\setminus B_{t}(\widetilde{x}_{0}) } \frac{|\wu(y)-(\wu)_{B_{t}(\widetilde{x}_{0})}|^q}{|y-\widetilde{x}_0|^{N+sq}} \thinspace dxdy\leq C
    \end{aligned}
\end{equation}
for every $s<\beta<1$ where we choose $\beta=(1+s)/2.$
Applying Lemma \ref{ccplem}, we get 
  \begin{equation*}
       \begin{aligned}
           \dashint_{B^{+}_{\rho/2}(\widetilde{x}_0)} |D\tu|^{p}\thinspace dx & +  \int_{B_{\rho/2}(\widetilde{x}_0)}  \dashint_{B^{+}_{\rho/2}(\widetilde{x}_0)} \frac{|\tu(x)-\tu(y)|^q}{|x-y|^{N+sq}} \thinspace dx \thinspace dy
        \\& \leq c\rho^{-p} \dashint_{B^{+}_{\rho}(\widetilde{x}_0)} |\tu|^{p} + c \int_{\RR^N\setminus B_{\rho}(\widetilde{x}_0)}\frac{|\tu(x)-\tu(y)|^{q}}{|x-y|^{N+sq}}
        + \dashint_{B^{+}_{\rho} (\widetilde{x}_0)} \eta^{m} \ck(x)\wf\wu .
       \end{aligned}
   \end{equation*}
 Now we estimate
 \begin{equation*}
     \begin{aligned}
         \rho^{-p} \dashint_{B^{+}_{\rho}(\widetilde{x}_0)} |\tu|^{p} \leq 2^{p-1} \relax[\wu]^{q}_{C^{\beta}(\RR^N)} \rho^{(\beta-1)p}.
     \end{aligned}
 \end{equation*}
Thus by choosing $\beta>1-\lambda>s$ we have \eqref{uleme3}.
\end{proof}
We will establish the H\"older continuity of the gradient of $\tu$ using a  perturbative method, specifically, by comparing it to the solution  $\th\in W_{\tu}^{1,p}(B^{+}_{\rho/4}(\widetilde{x}_0))$ of the Dirichlet problem given by:
\begin{equation}
    \label{theqn}
    \begin{aligned}
        \Div (\widetilde{A}(\widetilde{x}_0,D\th))=0 & \text{ in } B^{+}_{\rho/4}(\widetilde{x}_0),
        \\
        \th=\tu & \text{ on } \partial B^{+}_{\rho/4}(\widetilde{x}_0).
    \end{aligned}
\end{equation}
 It is important to note that the proof of \cite[Lemma 3.3]{AC23} does not rely on any specific properties of the data. Therefore, we can redo this theorem verbatim and obtain the following lemma.
\begin{lemma}\label{hlem}
    Let $\th$ be the solution of \eqref{theqn}. Then there exist constants $\sigma\in(0,1)$ and $C>0$ such that 
    \begin{equation}
        \begin{aligned}
                      &\dashint_{B^{+}_{\rho/4}(\widetilde{x}_0)} |D\th|^{p}\thinspace dx  \leq      \dashint_{B^{+}_{\rho/4}(\widetilde{x}_0)} |D\tu|^{p}\thinspace dx,
                \\
               & \|\th\|_{L^{\infty}(B^{+}_{\rho/4}(\widetilde{x}_0))} \leq       \|\tu\|_{L^{\infty}(B^{+}_{\rho/4}(\widetilde{x}_0))},
                \\
                &\Osc_{B^{+}_{\rho/4}(\widetilde{x}_0)} \th  \leq   \Osc_{B^{+}_{\rho/4}(\widetilde{x}_0)} \tu,
                \\
                 &\Osc_{B^{+}_{t}(\widetilde{x}_0)} D\th  \leq  C (\frac{t}{\rho})^{\sigma} \left\{\dashint_{B^{+}_{\rho/4}(\widetilde{x}_0)}|D\th|^{p}\thinspace dx 
 \right\}^{\frac{1}{p}}
 \text{ for all } t\in \left(0,\frac{\rho}{8}\right].
        \end{aligned}
    \end{equation}
\end{lemma}
Next, we will estimate all the integrals previously associated with the singular data to derive the boundary estimate for the gradient of $\tu - \th$. 
\begin{lemma}\label{dulem}
    Let us assume $\gamma<1.$  Let $\th$ and $\wu$ be the solution of \eqref{theqn} and \eqref{wueqn}, respectively. Then 
    \begin{equation}
        \dashint_{B^{+}_{\rho/4}(\widetilde{x}_0)} |D\tu - D\th|^{p}\thinspace dx \leq C \rho^{\Bar{\sigma}p} \text{ for some } \Bar{\sigma} \in (0,1).
    \end{equation}

\end{lemma}
\begin{proof}
Set $\tw=\tu-\th\in  W_{0}^{1,p}(B^{+}_{\rho/4}(\widetilde{x}_0)) $ and $\tw\equiv 0$ in $\RR^N\setminus B^{+}_{\rho/4}(\widetilde{x}_0).$ The new function $\tw\in \wsqq(\RR^N)$ which further implies that $\tw\in\wsq( B^{+}_{\rho/4}(\widetilde{x}_0))$ similar to \cite[ Lemma 2.4 and Lemma 5.2]{FM22}. By the definition of $\tw, $ along with the application of the last inequality of Lemma \ref{ulem} and the fact that $\dashint_{B^{+}_{\rho/4}(\widetilde{x}_0)} |D\th|^{p}\thinspace dx  \leq      \dashint_{B^{+}_{\rho/4}(\widetilde{x}_0)} |D\tu|^{p}\thinspace dx$ we get that 
 \begin{equation*}
     \begin{aligned}
         \dashint_{B^{+}_{\rho/2}(\widetilde{x}_0)} |D\tw|^{p}\thinspace dx \leq C_\lambda \rho^{-\lambda p} +\dashint_{B^{+}_{\rho} (\widetilde{x}_0)} \eta^{m} \ck(x)\wf\wu 
     \end{aligned}
 \end{equation*} 
 for any $\lambda>0.$  We set 
 \begin{equation*}
     \begin{aligned}
         V_{\mu}(z):=|z|^\frac{p-2}{2} z \text{ for } z\in\RR^n.
     \end{aligned}
 \end{equation*}
 Analogous to \cite[Equation (2.10)]{FM22} there exists a constant $C>0$ depending only on $n,p$ such that
 \begin{equation*}
     \begin{aligned}
         |  V_{\mu}(z_1)-  V_{\mu}(z_2)|^{2} \leq C (\widetilde{A}(\widetilde{x}_{0},z_1)-\widetilde{A}(\widetilde{x}_{0},z_2))\cdot(z_1-z_2) \text{ for all } z_1,z_2\in \RR^N.
     \end{aligned}
 \end{equation*}
Subsequently defining $\mathcal{V}^2:= |  V_{\mu}(D\widetilde{u})-  V_{D\widetilde{h}}(z_2)|^{2}$
we get that  
 \begin{equation*}
     \begin{aligned}
            \dashint_{B^{+}_{\rho/4}(\widetilde{x}_0)} \mathcal{V}^2 \thinspace dx\leq C(I_1+ I_2 + I_3 +I_4)
     \end{aligned}
 \end{equation*}
 where 
 \begin{equation*}
     \begin{aligned}
         I_1 &= \rho^{\alpha} \dashint_{B^{+}_{\rho/4}(\widetilde{x}_0)} |D\tu|^{p-1} |D \tw| \thinspace dx,
       \\ I_2 &= \dashint_{B^{+}_{\rho/4}(\widetilde{x}_0)} |\tf\tw|,
          \\I_3 &=\int_{B_{\rho/2}(\widetilde{x}_0)} \dashint_{B_{\rho/2}(\widetilde{x}_0)} \frac{|\tu(x)-\tu(y)|^{q-1}|\tw(x)-\tw(y)|}{|x-y|^{N+sq}} \thinspace dx \thinspace dy,
        \\  I_4 &=\int_{\RR^N \setminus B_{\rho/2}(\widetilde{x}_0)}\left( \dashint_{B_{\rho/2}(\widetilde{x}_0)} \frac{|\tu(x)-\tu(y)|^{q-1}|\tw(x)|}{|x-y|^{N+sq}} \thinspace dx \right)\thinspace dy
     \end{aligned}
 \end{equation*}
 similarly to \cite[pp. 13]{AC23}. 
 Using H\"older regularity and the last inequality of Lemma \ref{ulem} we obtain: 
 \begin{equation*}
     \begin{aligned}
         I_1 &\leq C \rho^{\alpha} (\dashint_{B^{+}_{\rho/4}(\widetilde{x}_0)} |D\tu|^{p} \thinspace dx)^{\frac{p-1}{p}}(\dashint_{B^{+}_{\rho/4}(\widetilde{x}_0)}|D \tw|^p \thinspace dx)^\frac{1}{p},
         \\&\leq C_\lambda \rho^{\alpha}\left( C_\lambda \rho^{-\lambda p} +\dashint_{B^{+}_{\rho} (\widetilde{x}_0)} \eta^{m} \ck(x)\wf\wu \right).
     \end{aligned}
 \end{equation*}
Next we estimate $I_3$ and $I_4$ similarly to \cite[pp. 14-15]{AC23}. Plugging all the estimate we get that 
 \begin{equation}\label{vest}
     \begin{aligned}
            \dashint_{B^{+}_{\rho/4}(\widetilde{x}_0)} \mathcal{V}^2 \thinspace dx\leq &C_\lambda \rho^{\alpha}\left( C_\lambda \rho^{-\lambda p} +\dashint_{B^{+}_{\rho} (\widetilde{x}_0)} \eta^{m} \ck(x)\wf\wu \right)
         + \dashint_{B^{+}_{\rho/4}(\widetilde{x}_0)} |\tf\tw|\thinspace dx+C_{\beta,\lambda}\rho^{(\beta-s)(q-1)-\lambda}
         \\
         & \hspace{9.5cm}+ \rho^{\beta-s}
     \end{aligned}
 \end{equation}
 for each $s<\beta<1$ and $(1-s)>\lambda>0.$ 
 \\
 First observe that $\frac{C_1}{(y_{N})^{\gamma}} \leq \wf(y) \leq \frac{C_2}{(y_{N})^{\gamma}}$ and $\wu \leq y_{N}^{\beta}$ due to $u\in C^{\beta}(\overline{\Omega}).$ Also, recalling  that $\tw=\tu-\th$ and using Hopf's Lemma, which implies that $\th\geq C y_{N},$ we infer that $|\tu-\th|\leq Cy_{N}^{\beta}$ for $\rho$ small enough. Substituting these into the equation we get:
 \begin{equation*}
     \begin{aligned}
         &\dashint_{B^{+}_{\rho} (\widetilde{x}_0)} \eta^{m} \ck(x)\wf\wu \leq \frac{C}{|B^{+}_{\rho} (\widetilde{x}_0)|} \int_{B^{+}_{\rho} (\widetilde{x}_0)} y_{N}^{\beta-\gamma} \leq C \varrho^{\beta-\gamma},
         \\
         &\dashint_{B^{+}_{\rho/4}(\widetilde{x}_0)} |\tf\tw|\thinspace dx \leq \frac{C}{|B^{+}_{\rho} (\widetilde{x}_0)|} \int_{B^{+}_{\rho} (\widetilde{x}_0)} y_{N}^{\beta-\gamma} \leq C \varrho^{\beta-\gamma}.
     \end{aligned}
 \end{equation*}

    Plugging these estimates into \eqref{vest},  we get the following for $\rho<<1$: 
    \begin{equation*}
        \begin{aligned}
            \dashint_{B^{+}_{\rho/4}(\widetilde{x}_0)} \mathcal{V}^2 \thinspace dx\leq &C_\lambda \rho^{\alpha}\left( C_\lambda \rho^{-\lambda p} +\dashint_{B^{+}_{\rho} (\widetilde{x}_0)} \eta^{m} \ck(x)\wf\wu \right)
          + \dashint_{B^{+}_{\rho/4}(\widetilde{x}_0)} |\tf\tw|\thinspace dx+C_{\beta,\lambda}\rho^{(\beta-s)(q-1)-\lambda} \\ &\hspace{9.5cm}+ \rho^{\beta-s}
            \\& \leq  C_\lambda \rho^{\alpha}\left(C_\lambda \rho^{-\lambda p}+ \rho^{\beta-\gamma}\right)+ C_\lambda \rho^{\beta-\gamma}
            +C_{\beta,\lambda}\rho^{(\beta-s)(q-1)-\lambda} + \rho^{\beta-s}
            \\& \leq C_\lambda \rho^{\alpha-\lambda p}  + C \rho^{\beta-\gamma} +C_{\beta,\lambda}\rho^{(\beta-s)(q-1)-\lambda} + \rho^{\beta-s}. 
        \end{aligned}
    \end{equation*}
    
 Since $\gamma<1$ we can choose 
     \begin{equation*}
         \begin{aligned}
             &\lambda = \min \{\frac{\alpha}{2p},\frac{(1-s)(q-1)}{4},(1-s)\},
             \\
             &\max \{\frac{1+s}{2}, \gamma\}<\beta <1 
         \end{aligned}
     \end{equation*}
     to obtain 
     \begin{equation*}
         \begin{aligned}
                 \dashint_{B^{+}_{\rho/4}(\widetilde{x}_0)} \mathcal{V}^2 \thinspace dx&\leq C \rho^{\sigma_0 p}
         \end{aligned}
     \end{equation*}
     with 
     \begin{equation*}
         \begin{aligned}
             \sigma_0 = \frac{1}{p} \min \{\frac{\alpha}{2},\frac{(1-s)(q-1)}{4},\frac{1-s}{2}\}.
         \end{aligned}
              \end{equation*}

Applying \cite[(2.9)]{FM22} and H\"older inequality, one can deduce that 
\begin{equation*}
    \frac{1}{C} \dashint_{B^{+}_{\rho/4}(\widetilde{x}_0)} |D\tu-D\th|^{p} \leq \left\{
    \begin{aligned}
        & \dashint_{B^{+}_{\rho/4}(\widetilde{x}_0)} \mathcal{V}^2 \text{ if } p\geq 2,\\
        & \left(\dashint_{B^{+}_{\rho/4}(\widetilde{x}_0)} \mathcal{V}^2 \right)^\frac{p}{2} \left(\dashint_{B^{+}_{\rho/4}(\widetilde{x}_0)}\Big(|D\tu|^{2}+ |D\th|^{2}\Big)^\frac{p}{2}\right)^{\frac{2-p}{2}} \text{ if } p\in(1,2).
    \end{aligned} \right.
\end{equation*}

Now for $p\geq 2$ we can choose $\Bar{\sigma}=\sigma_0.$ For $p\in(1,2)$ we note that, 
\begin{equation*}
    \dashint_{B^{+}_{\rho/4}(\widetilde{x}_0)} |D\th|^{p}\thinspace dx  \leq      \dashint_{B^{+}_{\rho/4}(\widetilde{x}_0)} |D\tu|^{p}\thinspace dx.
\end{equation*}
This observation combined with \eqref{uleme3}, leads to
\begin{equation*}
    \begin{aligned}
    \dashint_{B^{+}_{\rho/4}(\widetilde{x}_0)}\Big(|D\tu|^{2}+ |D\th|^{2}\Big)^{\frac{p}{2}} &\leq C
     \left(\rho^{-\lambda p}
 + \dashint_{B_{\rho}^{+}(x_{0})} \wf \wu\eta^{m} \ck(x)\right)
 \\
 &\leq C  \left(\rho^{-\lambda p}+\rho^{\beta-\gamma}\right) \leq C \rho^{-\lambda p}
    \end{aligned}
\end{equation*}
for small values of $\rho.$  Therefore for $p\in (1,2)$ we infer that 
\begin{equation*}
    \begin{aligned}
       \dashint_{B^{+}_{\rho/4}(\widetilde{x}_0)} |D\tu - D\th|^{p}\thinspace dx \leq   C (\rho^{\sigma_0 p})^{p/2} \rho^{\frac{-\lambda p(2-p)}{2}}.
    \end{aligned}
\end{equation*}
By choosing $\lambda=\frac{\sigma_0 p}{2(2-p)}$ and $\Bar{\sigma}=\frac{\sigma_0 p}{4}$ we conclude our desired result.
\end{proof}
As per the convention, we define the average of a function as follows:
Let $B\subset\RR^N$ be a measurable subset with respect to a
Borel (non-negative) measure $\lambda_{0}$ in $\RR^N$, with bounded positive measure $0<\lambda_0(B)<\infty$, and with $b : B \mapsto\RR^{k} , k \geq 1$, is a measurable map, we denote the average by
\begin{equation*}
    \begin{aligned}
        (b)_{B}:=\dashint_{B} b(x) d\lambda_{0}(x).
    \end{aligned}
\end{equation*}
Now
observe that, with Lemma \ref{dulem} already established, we can employ the ideas of \cite[pp.15-16]{AC23} to conclude the following result.
\begin{proposition}\label{gradcampest}
    Let $\wu$ be a solution of \eqref{wueqn}. Then there exists a radius $\rho_0\in(0,1)$ and constant $C>0,$ $\sigma_1\in(0,1)$ such that 
    \begin{equation*}
        \sup_{\widetilde{x}_0\in \Gamma_{r_0/2}} \dashint_{B^{+}_{\rho}(\widetilde{x}_0)} |D\tu - (D\tu)_{B^{+}_{\rho}(\widetilde{x}_0)}|^p \thinspace dx \leq C \rho^{\sigma_1 p} \text{ for every } \rho\in (0, \rho_0].
    \end{equation*} 
\end{proposition}
In the end, we establish the gradient H\"older regularity of the solution.
\\
\textit{\underline{Proof of Theorem \ref{gradientholderregularity} }}
\\[3mm]
The result follows by combining the interior Campanato estimate from \cite[Section 7]{FM22} with the boundary estimate in Proposition \ref{gradcampest}. Applying a standard covering argument along with Campanato's characterization of Hölder spaces, we conclude the proof.\hfill\qed

\section{Sobolev and H\"older Regularity For Singular Nonlinearity}\label{regsingnonlinear}
In this section we consider $u$ as a weak solution to the problem with  mixed local-nonlocal term and bounded data as given below:
  \begin{equation}\label{regeq}
      \begin{aligned}
           \fp u +  \fqs u &= \frac{\kg(x)}{u^{\delta}}+f(x)\text{ in } \Omega,
        \\
        u&>0 \text{ in } \Omega,\\
        u &=0 \text{ in } \RR^N \setminus \Omega
      \end{aligned}
  \end{equation}
  where $f\in L^{\infty}(\Omega), f\geq 0 ,$ and $\gamma,\delta>0.$  Furthermore we assume $\kg(x)\in L^{\infty}_{\loc}(\Omega)$ is a weight function such that 
  \begin{equation*}
    C_1 \dis^{-\gamma}(x)\leq \kg(x) \leq C_2 \dis^{-\gamma}(x) \text{ in } \Omega
\end{equation*}
for some $C_1, C_2>0.$

\par
First we set
 \begin{equation}\label{aleqn}
        0<\al<\left\{ 
        \begin{aligned}
          &\min \{ \frac{p-\gamma}{p-1+\delta},s\} \text{ if } p<q,
          \\
          & \min \{ \frac{p-qs}{p-q}, \frac{p-\gamma}{p-1+\delta},s\} \text{ if } q<p
        \end{aligned}
          \right.
      \end{equation}
      and fix this $\al$ throuhgout this section.
 To approximate $u,$
 we consider the following auxiliary problem for $\varepsilon>0:$
 \begin{eqnarray}\label{auxpb}
     \begin{array}{rll}
          \fp \uep +  \fqs \uep &= \kge(x) (\uep+\varepsilon)^{-\delta}+f(x),\quad \uep >0 \thickspace &\text{ in } \Omega,
        \\
        \uep &=0 &\text{ in } \RR^N \setminus \Omega
     \end{array}
 \end{eqnarray}
where 
\begin{equation*}
     \kge(x):=\left\{\begin{aligned}
       &(\kg(x)^\frac{-1}{\gamma}+\varepsilon^\frac{1}{\al})^{-\gamma} \text{ if } \kg(x)>0,
       \\
       & 0 \hspace{7.8em}\text{ otherwise
       }
    \end{aligned}\right.
\end{equation*} 
for $0\leq \gamma<p.$ Then we have 
\begin{equation*}
    \begin{aligned}
        C_1(d(x)+\varepsilon^\frac{1}{\al})^{-\gamma}\leq \kge(x) \leq C_1(d(x)+\varepsilon^\frac{1}{\al})^{-\gamma} \text{ in } \Omega.
    \end{aligned}
\end{equation*}
First, we prove an auxiliary theorem to get the lower behavior of $\uep$ independent of $\varepsilon>0.$
     \begin{theorem}\label{wtheta}
    For every $\Theta>0$ there exists a unique solution $w_{\Theta}\in \wo(\Omega)\cap C^{1,\alpha}(\Bar{\Omega})$ of the following equation 
    \begin{equation}\label{wthetaeqn}
        \begin{aligned}
            \fp w_{\Theta} + \fqs w_{\Theta} &= \Theta,\quad w_{\Theta}>0 \text{ in } \Omega,
            \\
            w_{\Theta}&=0 \text{ on } \Omega^{c}.
        \end{aligned}
    \end{equation}
    Moreover, $w_{\Theta}\rightarrow 0$ in $C^{1,\alpha}(\Bar{\Omega})$ as $\Theta\rightarrow 0.$
\end{theorem}
\begin{sketch}
    The existence of a solution in $\wo(\Omega)$ is established through the standard minimization technique, and uniqueness is ensured by the monotonicity of the operator. Furthermore, uniform $L^{\infty}$ bounds can be obtained through similar calculations as in \cite[Proposition 2.1]{FM22}. Thus we can apply \cite[Theorem 1.1]{AC23} to infer that $\|w_\Theta\|_{C^{1,\alpha}} \leq C$, independent of $\Theta$. By the compact embedding result, we have $w_\Theta \rightarrow w$ in $C^{1,\alpha^{\prime}}$. Given the uniqueness of the solution, it follows that $w = 0$, and we arrive at the desired result.
\end{sketch}
\subsection{Sobolev Regularity}
First we prove the existence and the uniqueness of weak solutions for the auxiliary problem.
\begin{lemma}\label{ueplemma}
    For each $\varepsilon>0$ there exists the unique solution $\uep \in \wo(\Omega)$ to the problem \eqref{auxpb}. Moreover, the sequence $\uep$ is decreasing in $\varepsilon$ and for each $D \subset\subset\Omega,$ there exists $C_{D}>0$ such that for all $\varepsilon>0,$ we have $C_{D}\leq u_1 \leq u_\varepsilon$ in $D.$
\end{lemma}
\begin{sketch}
 For fixed $\varepsilon>0$ and for each $v\in L^{\min\{p,q\}}(\Omega)$ we consider the following problem
    \begin{equation}\label{weqn}
        \begin{aligned}
            \fp w  + \fqs w= \kge(x)(|v|+\varepsilon)^{-\delta}+f, \quad w>0 \text{ in } \Omega,\quad w=0 \text{ on } \partial\Omega.
        \end{aligned}
    \end{equation}
    By standard minimization technique we get there exists a unique solution $w\in \wo(\Omega)$ of \eqref{weqn}. The energy functional $J:\wop(\Omega)\cap\wsq(\Omega) \longrightarrow \RR$ defined by
    \begin{equation*}
        J(w):= \frac{1}{p}\int_\Omega |\nabla w|^{p}+ \frac{1}{q}\int_\Omega [ w]^{q}-\int_\Omega \kge(x)(|v|+\varepsilon)^{-\delta}-\int_\Omega f|v|
    \end{equation*}
    is continuous, strictly convex and coercive.
    We  define the solution operator $S:L^{\min\{p,q\}}(\Omega)\longrightarrow L^{\min\{p,q\}}(\Omega)$ as $S(v)=w$ where $w$ is the unique solution of \eqref{weqn}. Then
    \begin{equation*}
        \begin{aligned}
            \|S(u)\|^{p}_{W_{0}^{1,p}}=C\|\nabla w\|^{p}_{L^{p}} \leq C \int_\Omega (|\nabla w|^{p}+ [w]_{\wsq}^{q}) \leq C \varepsilon^{-\delta} |\Omega|^{\frac{p-1}{p}} \|w\|_{L^{p}}. 
        \end{aligned}
    \end{equation*}
    Similarly we get 
    \begin{equation*}
        \begin{aligned}
            \|S(u)\|^{q}_{W_{0}^{s,q}} \leq C \int_\Omega (|\nabla w|^{p}+ [w]_{\wsq}^{q}) \leq C \varepsilon^{-\delta} |\Omega|^{\frac{q-1}{q}} \|w\|_{L^{q}}. 
        \end{aligned}
    \end{equation*}
    Hence using the Schauder fixed point, the boundary regularity of the regular problem \cite[Theorem 4]{FM22}, the strong comparison Principle \cite[Corollary 1.3]{AC23} and weak comparison principle for the regular problem as in \cite[Page 9]{GKS21} we obtain the result.
\end{sketch}
\begin{lemma}
The sequence $\{\uep\}$ is uniformly bounded in $\wo$ if $\gamma+(\frac{1}{q}-s)(1-\delta)<1$ and $0<\gamma+\delta<1.$
If $\gamma+\delta> 1,$ then
there exists $\kappa\geq 1$ such that $\uep^{\kappa}$ is uniformly bounded in $\wo(\Omega).$ 
\end{lemma}
\begin{proof}
    First we consider the case $\delta+\gamma<1$ and $\gamma+(\frac{1}{q}-s)(1-\delta)<1.$ Using the weak formulation of \eqref{auxpb} with $\uep$ as a test function we get
         \begin{equation*}
           \begin{aligned}
                 \|\uep\|_{W_{0}^{s,q}}^{q}\leq \int_{\Omega}|\nabla\uep|^{p} + A_{q}(\uep,\uep, \RR^{2N})&=\int_\Omega \frac{\kge\uep}{(\uep+\varepsilon)^{\delta}}+ \int_\Omega f \uep
             \\   &\leq C \int_\Omega \dis^{s(1-\delta)-\gamma} \left(\frac{\uep}{\dis^s}\right)^{1-\delta}+ C \|f\|_{L^{\infty}}  \|\uep\|_{W_{0}^{s,q}}
                \\& \leq C\Bigg[\left(\int_\Omega \dis^{\{-\gamma+s(1-\delta)\}\frac{q}{q-1+\delta}}\right)^\frac{q-1+\delta}{q}\left(\int_\Omega \left(\frac{\uep}{\dis^s}\right)^{q} \right)^{\frac{1-\delta}{q}} \\ &\hspace{5cm}+ \|\uep\|_{W_{0}^{s,q}}\Bigg]
                \\ &\leq C [ \|\uep\|_{W_{0}^{s,q}}^{1-\delta}+ \|\uep\|_{W_{0}^{s,q}}].
           \end{aligned}
         \end{equation*}
         From this we conclude $\|\uep\|_{W_{0}^{s,q}}<\infty$ independent of $\varepsilon.$
          We also obtain 
           \begin{equation*}
           \begin{aligned}
                 \|\uep\|_{W_{0}^{1,p}}^{p}\leq \int_{\Omega}|\nabla\uep|^{p} + A_{q}(\uep,\uep, \RR^{2N})&=\int_\Omega \frac{\kge\uep}{(\uep+\varepsilon)^{\delta}}+ \int_\Omega f \uep
            \\    &\leq C \int_\Omega \dis^{(1-\delta)-\gamma} \left(\frac{\uep}{\dis}\right)^{1-\delta} +C \|f\|_{L^{\infty}}  \|\uep\|_{W_{0}^{1,p}}
            \\& \leq C\Bigg[\left(\int_\Omega \dis^{\{-\gamma+(1-\delta)\}\frac{p}{p-1+\delta}}\right)^\frac{p-1+\delta}{p}\left(\int_\Omega \left(\frac{\uep}{\dis}\right)^{p} \right)^{\frac{1-\delta}{p}} \\
            &\hspace{5cm}+ \|\uep\|_{\wop} \Bigg ]    \\ &\leq C [ \|\uep\|_{W_{0}^{1,p}}^{1-\delta}
            + \|\uep\|_{W_{0}^{1,p}}]
           \end{aligned}
         \end{equation*}
         where we have used that $0<\gamma+\delta<1$ and consequently $\int_\Omega \dis^{\{-\gamma+(1-\delta)\}\frac{p}{p-1+\delta}}<\infty.$ Together with the uniform $W^{s,q}$ bound, this implies $\{u_\epsilon\}$ is uniformly bounded in $\wo(\Omega).$
         \par
         Next we consider $\gamma+\delta>1$ and take $\uep^{\rho}$ for $\rho\geq  1$ ( will be chosen later) as the test  function in \eqref{auxpb}. Then we get 
 \begin{equation*}
    \begin{aligned}
        \int_{\Omega} |\nabla \uep|^{p-2} \nabla \uep \nabla \uep^{\rho} &+ \int_{\RR^{2N}}\frac{|\uep(x)-\uep(y)|^{q-2}(\uep(x)-\uep(y))(\uep^{\rho}(x)-\uep^{\rho}(y))}{|x-y|^{N+sq}} \thinspace dx\thinspace dy 
        \\
       & \leq  \int_\Omega \left(\frac{\kge(x)}{{(\uep+\varepsilon)}^{\delta}}+f(x)\right)\uep^{\rho}.
    \end{aligned}
\end{equation*}
Using \cite[Lemma A.2]{BP16} we get that $A_{q}(u, \uep^{\rho},\RR^{2N})\geq 0.$ Observe that 
\begin{equation*}
    \begin{aligned}
      \int_{\Omega} |\nabla \uep|^{p-2} \nabla \uep \nabla \uep^{\rho} =\frac{p\rho}{p+\rho-1} \int_{\Omega} |\nabla \uep^{\frac{p+\rho-1}{p}}|^{p}.
    \end{aligned}
\end{equation*}
Now 
\begin{equation*}
    \begin{aligned}
        \int_\Omega f(x)\uep^{\rho}\leq C \|\uep^{\frac{p+\rho-1}{p}}\|_{\wop}^{\frac{\rho p}{p+\rho-1}}
    \end{aligned}
\end{equation*}
and choose $\rho\geq \max\left\{ 1,\frac{(p-1)(\gamma-1+\delta)}{p-\gamma}\right\}$ to get 
\begin{equation*}
    \begin{aligned}
         \int_\Omega \frac{\kge(x)}{{(\uep+\varepsilon)}^{\delta}}\uep^{\rho}\leq \int_\Omega \dis^{-\gamma} \uep^{\rho-\delta} \leq \int_{\Omega} \dis^{-\gamma+\frac{(\rho-\delta)p}{p+\rho-1}} \left(\frac{\uep^{\frac{p+\rho-1}{p}}}{\dis}\right)^{\frac{(\rho-\delta)p}{p+\rho-1}}\leq C \|\uep^{\frac{p+\rho-1}{p}}\|_{\wop}^{\frac{(p+\rho-1)p}{\rho-\delta}}.
    \end{aligned}
\end{equation*}
Hence $\uep^{\frac{p+\rho-1}{p}}$ is bounded in $\wopo.$
Next for $\varrho\geq 1$ (chosen later) we use $\uep^{q(\varrho-1)+1}$ as the test function and observe that 
\begin{equation*}
    \begin{aligned}
         \int_{\Omega} |\nabla \uep|^{p-2} \nabla \uep \nabla \uep^{q(\varrho-1)+1} =[q(\varrho-1)+1] \int_{\Omega} |\nabla \uep|^{p} \uep^{q(\varrho-1)} \geq 0.    \end{aligned}
\end{equation*}
Hence we get that 
\begin{equation*}
    \begin{aligned}
        \|\uep^{\varrho}\|_{\wsq}^{q}\leq \int_\Omega \left(\frac{\kge(x)}{{(\uep+\varepsilon)}^{\delta}}+f(x)\right)\uep^{q(\varrho-1)+1}.
    \end{aligned}
\end{equation*}
Obviously
\begin{equation*}
    \begin{aligned}
        \int_{\Omega} f(x)\uep^{q(\varrho-1)+1} \leq \| \uep^{\varrho}\|_{L^{q}}^{\frac{\varrho q}{q(\varrho-1)+1}}
    \end{aligned}
\end{equation*}
and choose $\varrho\geq \max\left\{1, \frac{(s-\frac{1}{q})(q+\delta-1)}{sq-\gamma}\right\}$ to get 
\begin{equation*}
    \begin{aligned}
          \int_\Omega \frac{\kge(x)}{{(\uep+\varepsilon)}^{\delta}}\uep^{q(\varrho-1)+1}&\leq \int_{\Omega} \dis^{-\gamma+s\frac{q(\varrho-1)+1-\delta}{\varrho}} \left(\frac{\uep^{\varrho}}{\dis^{s}}\right)^{\frac{q(\varrho-1)+1-\delta}{\varrho}}
          \\
          &\leq C \|\uep^{\varrho}\|_{\wsqo}^{\frac{q(\varrho-1)+1-\delta}{\varrho}}. 
    \end{aligned}
\end{equation*}
Thus $\uep^{\varrho}$ is bounded in $\wsqo.$
Combining we choose 
\begin{equation*}
    \begin{aligned}
        \kappa\geq \max\left\{1, (p-1)(\frac{1}{p}+\frac{\gamma-1+p}{p-\gamma}), \frac{(s-\frac{1}{q})(q+\delta-1)}{sq-\gamma} \right\}
    \end{aligned}
\end{equation*}to get that $\uep^{\kappa}\in\wo$ and conclude the result.
\end{proof}
    \begin{corollary}\label{existenceofminimialsol}
    Up to a subsequence $\uep$ pointwise converges to $u,$ where $u$ is the minimal solution of \eqref{regeq}.
\end{corollary}

\subsection{Gradient H\"older Boundary Regularity}
In this section we derive the boundary gradient H\"older regularity of $u$ where $u$ is a weak solution  to \eqref{regeq} in the sense of Definition \ref{soldef} for $\gamma+(\frac{1}{q}-s)(1-\delta)<1$ and $0<\gamma+\delta<1.$ 
\\
\textit{\underline{Proof of Theorem \ref{gradhld}}}
\\[3mm]
Let $w_\Theta \in \wopo \cap \wsqo$ be a solution to \eqref{wthetaeqn} for some sufficiently small $\Theta(\Omega, \gamma, \delta) > 0$. Then, $w_\Theta$ serves as a subsolution to \eqref{regeq}. Note that $u_1 \in \wopo \cap \wsqo$, where $u_1$ solves \eqref{auxpb} for $\varepsilon = 1$. By the weak comparison principle, it follows that $w_\Theta \leq u_1$ in $\Omega$. Consequently, by Hopf's Lemma \cite[Theorem 1.2]{AC23} and the monotonicity of $\{ \uep \}_\varepsilon$, we deduce  
\begin{equation*}
    C_1 \dis(x) \leq u_1 \leq \uep \quad \text{in } \Omega.
\end{equation*}
Define $\dis_e(x)$ as  
\begin{equation*}
\dis_e(x) = 
\begin{cases} 
\dis(x) & \text{if } x \in \Omega, \\
-\dis(x) & \text{if } x \in (\Omega^c)_{\varepsilon^{1/\alpha}}, \\
-\varepsilon^{1/\alpha} & \text{otherwise} 
\end{cases}
\end{equation*}
where $\varepsilon \in (0, \varepsilon^*)$. Additionally, define  
\begin{equation*}
\Bar{\omega_\rho}(x) = 
\begin{cases}
\big( \dis_e(x) + \varepsilon^{\frac{1}{\alpha}} \big)_+^\alpha & \text{in } \Omega \cup (\Omega^c)_{\varepsilon^{1/\alpha}}, \\
0 & \text{otherwise}
\end{cases}
\end{equation*}
for $\frac{qs}{q-1} < \alpha < 1$. It follows that  
\begin{equation*}
    \begin{aligned}
     &\fp \Bar{\omega_\rho} \geq C \big( \dis(x) + \varepsilon^{\frac{1}{\alpha}} \big)^{\alpha p - \alpha - p} \quad \text{in } \Omega_\rho, \text{ and }   \\[1.5em]
     &\fqs \Bar{\omega_\rho} = h \quad \text{in } \Omega_\rho, \quad \text{where } h \in L^\infty(\Omega_\rho), \text{ independent of } \varepsilon.
    \end{aligned}
\end{equation*}
For sufficiently small $\rho$ and $\varepsilon^*$, we deduce 
\begin{equation*}
   \begin{aligned}
        \fp (\Gamma \wro) + \fqs (\Gamma \wro) &\geq C \Gamma^{p-1} (\dis(x)+ \varepsilon^{\frac{1}{\alpha}})^{\alpha p - \alpha -p} -\Gamma^{q-1} \|h\|_{L^{\infty}(\Omega_\rho)}\\
    &\geq \Gamma^{p-1} \frac{C}{2} (\dis(x)+ \varepsilon^{\frac{1}{\alpha}})^{\alpha p - \alpha -p}.
   \end{aligned}
\end{equation*}

Observe that $\frac{qs}{q-1} < \alpha < 1 < \frac{p-\gamma}{p-1+\delta}.$ Thus there exists $r_0 > 0$ such that  
\begin{equation*}
\alpha p - \alpha - p \leq -\gamma - \alpha \delta - r_0 < -\gamma - \alpha \delta.
\end{equation*}
Therefore, for sufficiently small $\rho$, and $\dis(x) \ll 1$ 
\begin{equation}
    \begin{aligned}
         \fp (\Gamma \wro) + \fqs (\Gamma \wro) &\geq \Gamma^{p-1} \frac{C}{2} (\dis(x)+ \varepsilon^{\frac{1}{\alpha}})^{\alpha p - \alpha -p}
         \\&\geq \Gamma^{p-1+\delta} \frac{C}{2}(\Gamma(\dis+\varepsilon^{\frac{1}{\alpha}})^{\alpha})^{-\delta}\kge(x)(\dis+\varepsilon^\frac{1}{\alpha})^{-r_{0}}
    \end{aligned}
\end{equation}
in $\Omega_\rho.$ Thus
\begin{equation*}
    \begin{aligned}
       \kge(x) \Gamma^{p-1+\delta}\frac{C}{2}(\Gamma(\dis+\varepsilon^{\frac{1}{\alpha}})^{\alpha})^{-\delta}(\dis+\varepsilon^\frac{1}{\alpha})^{-r_{0}} &\geq \kge(x)  \Gamma^{p-1+\delta}\frac{C}{2}(\Gamma(\dis+\varepsilon^{\frac{1}{\alpha}})^{\alpha})^{-\delta}+\|f\|_{L^{\infty}}
        \\&\geq \kge(x)(\Gamma(\dis+\varepsilon^\frac{1}{\alpha})^\alpha+\varepsilon)^{-\delta} + f
        \\&\geq\kge(x)(\Gamma\wro+\varepsilon)^{-\delta} +f 
    \end{aligned}
\end{equation*}
choosing $\rho^{*},$ $\varepsilon^{*}$ small enough and  $\Gamma > 1$ large enough such that $\Gamma^{p-1+\delta} \frac{C}{2} \geq 1.$  Now using $\uep\geq C \dis(x)$ in $\Omega,$ we have that $\frac{\kge}{(\uep+\varepsilon)^{\delta}}\leq C\kg \dis^{-\delta}$ for some $C$ independent of $\varepsilon.$ Since $C\kg \dis^{-\delta}\in L^{\infty}_{\loc}(\Omega)$ we can use Theorem \ref{interiorregularity} to get that $\|\uep\|_{L^{\infty}(\Omega\setminus\Omega_{\rho})}\leq C_{\rho}$ and the bound is independent of $\varepsilon.$ Consequently, choosing $\Gamma > 1$ sufficiently large, independent of $\varepsilon$, ensures  
\begin{equation*}
    \begin{aligned}
        \uep\leq \|\uep\|_{L^{\infty}(\Omega\setminus\Omega_{\rho})}\leq C_{\rho} \leq \Gamma(\dis+\varepsilon^\frac{1}{\alpha})^{\alpha}=\Gamma\wro \text{ in } \Omega\setminus\Omega_{\rho}.
    \end{aligned}
\end{equation*}
By the weak comparison principle \cite[Proposition 4.1]{AC23}, we deduce $\uep \leq \Gamma \Bar{\omega_\rho}$ in $\Omega$ for all $\varepsilon \in (0, 1)$. Finally, passing to the limit as $\varepsilon \to 0$ yields that $u(x)\leq \Gamma\dis^{\alpha}(x)$ in $\Omega$ for any $\alpha\in(\frac{qs}{q-1},1)$. Now since $u>0$ we can find $\Theta(\Omega,\gamma,\delta)>0$ small enough such that $0<\Theta\leq \frac{\kg}{u^{\delta}}+f(x).$ Now for that $\Theta$ let $w_{\Theta}$ solves \eqref{wthetaeqn}. From Hopf's Lemma \cite[Theorem 1.2]{AC23} and Theorem \ref{wcp} we obtain $ C \dis(x)\leq u$ in $\Omega.$ Thus, $0<\Theta\leq \frac{\kg}{u^{\delta}}+f(x)\leq \frac{C}{\dis^{\gamma+\delta}}$ in ${\Omega}$ for $C=C(\Omega,\|f\|_{L^{\infty}})$ large enough. Now from  Theorem \ref{gradientholderregularity} we get the desired result.\hfill\qed

\subsection{Boundary H\"older Regularity For Strongly Singular Problem}
Now we focus on the strongly singular nonlinearity, i.e. the case $\gamma+\delta>1.$ 
\\
\textit{\underline{Proof of Theorem \ref{strongsingreg}}}
\\[3mm]
Let $w_\Theta\in \wopo\cap\wsqo$ is a solution to \eqref{wthetaeqn}  for some $\Theta(\Omega,\gamma,\delta)>0$ sufficiently small. Then $w_{\Theta}$ is a subsolution to \eqref{regeq}. Observe that $u_1\in \wopo\cap\wsqo$ where $u_1$ solves \eqref{auxpb} for $\varepsilon=1.$ Therefore by the weak comparison principle we get that $w_\Theta\leq u_1$ in $\Omega.$ Consequently by the Hopf's Lemma \cite[Theorem 1.2]{AC23} and the fact that  $\{\uep\}_\varepsilon$ is decreasing, we obtain 
\begin{equation*}
    C_1 \dis(x) \leq u_1 \leq \uep \text{ in } \Omega.
\end{equation*}

 Define $\dis_e(x)$  as
\begin{equation*}
    \dis_{e}(x)= \left\{
    \begin{aligned}
        \dis(x) &\text{ if } x \in \Omega,\\
        -\dis(x) &\text{ if } x \in (\Omega^{c})_{\rho},\\
        -\rho &\text{ otherwise }
    \end{aligned}
    \right. 
\end{equation*}
and
 \begin{equation*}
\Bar{\omega}_\rho(x)=\left\{
\begin{aligned}
    &\Gamma(\dis_{e}(x)+ \va^{\frac{1}{\al}})_{+}^{\al} \text{ in } \Omega \cup (\Omega^{c})_{\rho},
    \\
    &0 \text{ otherwise}
\end{aligned}\right.    
\end{equation*}     
   for some $\Gamma>0$ (chosen later), $\rho\in(0,\rho_1)$ and $\al$ is defined in \eqref{aleqn}. Observe that $\al p - \al -p\leq -qs+\al(q-1) .$
      Then using \cite[Theorem 3.11]{GKS23}, there exist $\va_1,\rho_1>0$ such that for all $\va \in (0,\va_1)$ and $\rho \in (0, \rho_1),$   
   \begin{equation*}
\begin{aligned}
        \fp \Bar{\omega}_\rho + \fqs  \Bar{\omega}_\rho \geq  C \al^{p-1}(1-\al)(p-1)(\dis(x)&+ \varepsilon^{\frac{1}{\al}})^{\al p - \al -p} + C \Gamma^{q-1} (\dis(x)+ \varepsilon^{\frac{1}{\al}})^{-qs+\al(q-1)} 
    \\
    & \geq C (\dis(x)+ \varepsilon^{\frac{1}{\al}})^{\al p - \al -p} 
    \text{ near the boundary. }
\end{aligned}
\end{equation*}
Since $\gamma+\delta>1$, we have $\frac{p-\gamma}{p-1+\delta}<1.$ Since $\al<\frac{p-\gamma}{p-1+\delta}$, we get $\al p-\al-p< -\gamma-\al\delta\leq 0.$ Hence,
\begin{equation*}
    \begin{aligned}
C  (\dis(x)+ \varepsilon^{\frac{1}{\al}})^{\al p - \al -p}    &\geq C  \kge(x) (\dis+\varepsilon^{\frac{1}{\al}})^{-\al \delta} \geq (\Gamma(\dis+\varepsilon^\frac{1}{\al})^\al)^{-\delta} \\
  &\geq \kge(x)(\Gamma(\dis+\varepsilon^\frac{1}{\al})^\al+\varepsilon)^{-\delta}=\kge(x)(\Gamma\wro+\varepsilon)^{-\delta}
    \end{aligned}
\end{equation*}
where we choose $\Gamma>1$ large enough such that  $\Gamma^{\delta}  C \geq 1.$ Now using $\uep\geq C \dis(x)$ in $\Omega,$ we have that $\frac{\kge}{(\uep+\varepsilon)^{\delta}}\leq C\kg \dis^{-\delta}$ for some $C$ independent of $\varepsilon.$ Since $C\kg \dis^{-\delta}\in L^{\infty}_{\loc}(\Omega)$ we can use Theorem \ref{interiorregularity} to get that $\|\uep\|_{L^{\infty}(\Omega\setminus\Omega_{\rho})}\leq K_{\rho}$ and the bound is independent of $\varepsilon.$ Consequently we chose $\Gamma>1$ large enough independent of $\varepsilon$ to get that 
\begin{equation*}
    \begin{aligned}
        \uep\leq \|\uep\|_{L^{\infty}(\Omega\setminus\Omega_{\rho})}\leq K_{\rho} \leq \Gamma(\dis+\varepsilon^\frac{1}{\al})^{\al}=\Gamma\wro \text{ in } \Omega\setminus\Omega_{\rho}.
    \end{aligned}
\end{equation*}
Then using \cite[Proposition 4.1]{AC23} we get that $\uep\leq \Gamma\wro$ in $\Omega$ for all $\varepsilon\in(0,1). $ Passing to the limit as $\varepsilon\rightarrow 0$, we deduce that $u\leq C \dis^{\al}(x)$ in $\Omega.$ Since $u\geq 0$ is bounded and $u^{\kappa}\in \wo(\Omega)$ then we have 
\begin{equation*}
    \begin{aligned}
        \int_{\RR^{N}} \frac{u^{q}}{(1+|x|)^{N+sq}}\lesssim \int_{\RR^{N}} \left(\frac{u^{\kappa q}}{(1+|x|)^{N+sq}}\right)^{1/\kappa}<\infty.
    \end{aligned}
\end{equation*}
Hence $u\in \mathcal{X}(\Omega)$ then we can combine this boundary behaviour and the interior H\"older regularity (Theorem \ref{interiorregularity}) to  conclude the result. \hfill\qed

\begin{remark}
    The above calculation demonstrates that if  $0 < \gamma < p,$ $\gamma + \delta > 1$ and $u$ is a minimal solution to \eqref{regeq} with $f(x) \equiv 0$ , then we have $u \in C^{\al}(\overline{\Omega})$.\\ 
    Furthermore, for any $0<\gamma<p,$ a weak minimal solution $u$ to \eqref{regeq}  belongs to the solution space  $\mathcal{X}(\Omega) \cap \mathcal{C}_{\dis_{\gamma,\delta}}.$ The conical shell $\mathcal{C}_{\dis_{\gamma,\delta}}$ is defined as the set of all continuous functions $u$ in $\Omega$ such that
    \begin{itemize}
        \item[(a)] $C_1 \dis(x) \leq u(x) \leq C_2 \dis^{\xi}(x)$ for any $\xi\in(\frac{qs}{q-1},1)$ if $\gamma+\delta<1$ and $\gamma+(\frac{1}{q}-s)(1-\delta)<1,$
        \item[(b)] $C_3 \dis(x) \leq u(x) \leq C_4 \dis^{\al}(x)$ for $\gamma+\delta>1.$
    \end{itemize}
    In the light of remark \ref{cndgmdel} the condition $\gamma+(\frac{1}{q}-s)(1-\delta)<1$ can be dropped if $q<p.$

\end{remark}

\section{Comparison Principles}\label{comp}
In this section, we aim to infer weak and strong comparison principles for the mixed local-nonlocal operator with singular term and nonnegative data. Examining the proofs of \cite[Theorem 1.5]{GKS21} and \cite[Proposition 4.3]{GKS23}, the proof of the weak comparison principle follows verbatim.
\begin{theorem}[Weak Comparison Principle]\label{wcp}
  Let us assume $\gamma\leq \min\{1+s-\frac{1}{q}, 2-\frac{1}{p}\}.$ Then if $\underline{w},\overline{w}\in \mathcal{X}(\Omega)$ be sub and supersolution of \eqref{regeq} respectively then $\underline{w}\leq \overline{w}$ a.e. in $\Omega.$ 
\end{theorem}
As a consequence of the Theorem \ref{wcp} we get the following uniqueness result.
\begin{corollary}
    If $\gamma\leq \min\{1+s-\frac{1}{q}, 2-\frac{1}{p}\},$ then the equation \eqref{regeq} has at most one solution.
\end{corollary}
 To prove the strong comparison principle we need some estimates for $p$ Laplacian and fractional $q$-Laplacian which will be discussed in the next two lemmas.
\begin{lemma}\label{plapest}
    Let $f\in C_{c}^{2}(D)$ where $D$ is any  open bounded set in $\Omega$ and $a\in[-1,1].$  Then, for any given $v\in \wo(\Omega)\cap C^{1,\alpha}(\Bar{\Omega}),$ there exist $C$ independent of $a$ such that 
    \begin{equation*}
      \left |(\fp) (v-af)-(\fp) v\right|\leq \left\{ \begin{aligned}
            &C \,|a|^{p-1}\text{ for } 1<p\leq 2,\\
            & C \, |a| \hspace{1.4em}\text{ for } p\geq 2.
        \end{aligned}\right.
    \end{equation*}
   in the weak sense on the $\supp\, f.$ 
\end{lemma}
\begin{proof}
    Fix $f\in C_{c}^{2}(D)$ and denote $U:=\supp \,f.$  First we consider the case: $1<p\leq 2$ and from \cite[Lemma 2]{Lin16} recall that, for all $t\in(0,1]$ there is $C>0$ such that 
    \begin{equation*}
        \begin{aligned}
            |[m+n]^{t}-[n]^{t}| \leq C|n|^{t} \text{ for all } m,n\in \RR.
        \end{aligned}
    \end{equation*}
    For any $\varphi\in \wo(U), \ \varphi\geq 0$ we have 
    \begin{equation*}
        \begin{aligned}
            \left|\int_{\Omega} ([\nabla (v-af)]^{p-1}-[\nabla v]^{p-1}) \nabla \varphi\right|&\leq   \int_{\Omega} \left|([\nabla (v-af)]^{p-1}-[\nabla v]^{p-1})\right| \left|\nabla \varphi\right|
            \\&\leq \int_{\Omega} C|a|^{p-1}|\nabla f|^{p-1} |\nabla\varphi|.
        \end{aligned}
    \end{equation*}
   Since $f\in C_{c}^{2}(\Omega)$ we can now conclude this case using the H\"older inequality.
 Next consider $p\geq 2$ and from \cite[Lemma 2.1]{Jar18} recall that 
   \begin{equation}\label{pgetwoalgine}
       \begin{aligned}
           (m+n)^{t}\leq \max\{1, 2^{t-1}\}(m^{t}+n^{t}) \text{ if } m,n\geq 0, t>0.
       \end{aligned}
   \end{equation}

Moreover we infer 
   \begin{equation*}
       \begin{aligned}
           \left|\int_{\Omega} ([\nabla (v-af)]^{p-1}-[\nabla v]^{p-1}) \nabla \varphi\right|&\leq |a|\int_{\Omega}(p-1)|\nabla \varphi||\nabla f|\int_{0}^{1} |\nabla v - a t \nabla f|^{p-2}\thinspace dt
           \\ & \leq C|a|\int_{\Omega} |\nabla\varphi||\nabla f|\int_{0}^{1} \left(|\nabla v|+t|a||\nabla f|\right)^{p-2} \thinspace dt
           \\
           & \leq C |a| \int_{\Omega} |\nabla \varphi||\nabla f| (|\nabla v|^{p-2}+|\nabla f|^{p-2})
       \end{aligned}
   \end{equation*}
where we have used \eqref{pgetwoalgine} in the last line. We conclude the result using the fact that $v\in \wo(\Omega)\cap C^{1,\alpha}(\Bar{\Omega})$ and $f\in C_{c}^{2}(\Omega).$
\end{proof}
\begin{lemma}\label{estmix}
    Let \( D \subset \mathbb{R}^N \) be an open bounded set, let \( K \subset\subset \mathbb{R}^N \setminus D \) with \( |K| > 0 \), \( \eta \in (0, 1] \), \( s \in (0, 1) \), \( p > 1 \) and $q>\frac{1}{1-s}$. 
Fix \( f \in C_c^2(D) \) with \( 0 \leq f \leq 1 \). For any $v\in \wo(\Omega)\cap C^{1,\alpha}(\Bar{\Omega}),$ there exists  $a_0,b>0$ such that for all $a\in (0,a_0)$
\begin{equation*}
    \begin{aligned}
        \fp (v- a f -\eta 1_{K}) + \fqs  (v- a f -\eta 1_{K}) \geq \fp v + \fqs v + b \text{ in } \supp f.
    \end{aligned}
\end{equation*}
\end{lemma}
\begin{proof}
   Since $K \cap \supp f =\varnothing$ , from Lemma \ref{plapest} for all $a\in(0,1]$  we get
   \begin{equation*}
       \begin{aligned}
            \fp (v- a f -\eta 1_{K}) =  \fp (v- a f) \geq \fp v - C \max\{a, a^{p-1}\}  \mbox{ in } \supp f
       \end{aligned}
   \end{equation*}
   for some $C>0.$ Using  \cite[Lemma 3.4 and Lemma 3.5]{Jar18},  for all $a\in(0,1],$ we have
    \begin{equation*}
       \begin{aligned}
            \fqs (v- a f -\eta 1_{K}) \geq \fqs v -  C_1\max\{a, a^{q-1}\} + \widetilde{C} \min\{\eta, \eta^{q-1}\}
       \end{aligned}
   \end{equation*}
   for some $C_1, \widetilde{C}>0.$ Now we fix $a_{0}(\eta)\in(0,1]$ such that 
   \begin{equation*}
       \begin{aligned}
           b= - \max\{a, a^{p-1}\} C - \max\{a, a^{q-1}\} C_1 + \widetilde{C} \min\{\eta, \eta^{q-1}\}>0
       \end{aligned}
   \end{equation*}
   to conclude the result.
\end{proof}
We now present the proof of the strong comparison principle for the mixed local-nonlocal operator with a singular term.
\\
\textit{\underline{Proof of Theorem \ref{scp} }}
\\[3mm]
 Using Lemma \ref{estmix}, we adopt a proof strategy similar to that of \cite[Theorem 2.19]{GKS23}. However, for the sake of completeness, we provide a sketch of the proof here. \\
Let $u,v$ be $C^{1,\alpha}(\overline{\Omega})$ solutions of \eqref{uvsol}. By weak comparison principle, we know that $v\geq u$ in $\Omega.$ Let $V$ be a non-empty open subset of $\Omega$ where $\inf_{V} (v-u) >0.$ Clearly such a $V$ exists, otherwise this implies $v\equiv u$ in $\Omega.$ Define $\eta:=\inf_{V}(v-u).$ Without loss of generality, assume that $0<\eta\leq 1.$ Fix $E\subset\subset\Omega\setminus V$ and $f\in C_{c}^{2}(\Omega\setminus V)$ such that $f\equiv 1$ in $E$ and $0 \leq f \leq 1.$ Define $a_{0}, b>0$ as given by Lemma \ref{estmix}. For all $a\in(0,a_{0}]$ with $a< L\, \dis(\supp\, f, \partial\Omega),$ set $\omega=(v- a f -\eta 1_{V}).$ Since $v\in \wo(\Omega)\cap C^{1,\alpha}(\Bar{\Omega}),$ by Lemma \ref{estmix} we infer 
   \begin{equation*}
       \begin{aligned}
            \fp \omega + \fqs  \omega &\geq \fp v + \fqs v + b
            \\ & \geq  \lambda\omega^{-\delta} + g +b + \lambda(v^{-\delta}- \omega^{-\delta} )
            \\ &\geq \lambda\omega^{-\delta} + g +b + \lambda(v^{-\delta}- (v- a)^{-\delta} ) \mbox{ in } \supp f.
       \end{aligned}   \end{equation*}
        Now we choose $a>0$ sufficiently small such that $b + (v^{-\delta}- (v- a)^{-\delta} )>0$ in $\supp f.$ This gives
       \begin{equation*}
           \begin{aligned}
              \fp \omega &+ \fqs  \omega \geq  \lambda\omega^{-\delta} + g > \lambda\omega^{-\delta} + f  \text{ in } \supp f,
           \end{aligned}
       \end{equation*}
       and $\omega\geq u$ in $\RR^{N}\setminus\supp f.$ Hence using the weak comparison principle we get that $v\geq u + a$ in $E.$ Thus,  $\inf_{E}(v-u)>0$ for all $E \subset\subset\Omega\setminus V$. Since $E$ is arbitrary we conclude our result.
\hfill $\square$
 
\begin{remark}
   To the best of our knowledge, the result concerning the strong comparison principle is novel, even for the homogenous mixed local-nonlocal operator (i.e., \( p = q  \)) with a singular nonlinearity.
    \end{remark}

\section{Sublinear Perturbation Problem}\label{picone}
 In this section, we consider the problem
\begin{equation}\label{sublineareqn}
    \begin{aligned}
        \fp u + \fqs u &=\frac{\lambda}{u^{\delta}}+u^{l} \text{ in }\Omega,
        \\
        u&=0 \text{ in } \Omega^c
        \end{aligned}
\end{equation}
where $0<l<\min\{p-1,q-1\}$ and $0<\delta<1.$ 
Our aim is to prove the global boundedness of the weak solutions.
\begin{lemma}\label{uobd}
    Let $\uo\in\wo(\Omega)$ be a solution to the problem
    \begin{equation}\label{uoeqn}
        \left\{
        \begin{array}{rll}
            \fp \uo +  \fqs \uo &= F(x,\uo);\;\; \uo>0 &\text{ in } \Omega,
         \\   \uo&=0 &\text{ in } \Omega^{c}
        \end{array}\right.
    \end{equation}
    where $|F(x,s)|\leq \Lambda(s^{-\delta}+s^{l})$ for $l\leq \min\{p-1,q-1\}.$  Then $\|\uo\|_{L^{\infty}(\Omega)}$ is bounded depending on $\Lambda, l,N,p$ and $\|\uo\|_{L^{l}(\Omega)}.$  
\end{lemma}

\begin{proof}    
Let $\alpha<1$ be a constant to be chosen later and define $\wuo=\alpha\uo.$ Define $v_k:=\wuo-c_k$ for $c_k:=1-2^{-k},$  $w_k:=(v_k)_{+}$ and $U_k:=\|w_k\|_{L^{l}}^{l}.$ Clearly, $w_k \in \wopo\cap\wsqo.$ Now we can use $w_k$ as our test function in the weak formulation of \eqref{uoeqn} and get 
\begin{equation}\label{testfn}
    \begin{aligned}
        \int_\Omega |\nabla\wuo|^{p-2} \nabla\wuo \nabla w_k &+ A_{q}(\wuo,w_k,\RR^{2N})  \leq \left[\int_\Omega C(\wuo^{-\delta}+ \wuo^{l})w_k\right]
    \end{aligned}
\end{equation}
since $\wuo\leq \uo,$  $\alpha<1$ and $\min\{p-1-l,q-1-l\}\geq 0.$ 
From \cite[Equation (14)]{FP14} we get that, for any measurable function $z$ and for almost every pair of points $x,y\in\RR^N,$ we have
\begin{equation*}
    |z_{+}(x)-z_{+}(y)|^{q}\leq |z(x)-z(y)|^{q-1}(z_{+}(x)-z_{+}(y)).
\end{equation*}  
Hence,
\begin{equation} \label{sqnorm}
    \begin{aligned}
        |(w_k(x)-w_k(y))|^{q}\leq |\wuo(x)-\wuo(y)|^{q-2}(\wuo(x)-\wuo(y))(w_k(x)-w_k(y)).
    \end{aligned}
\end{equation}
Moreover 
\begin{equation} \label{pnorm}
    \begin{aligned}
         \int_\Omega |\nabla\wuo|^{p-2} \nabla\wuo \nabla w_k = \int_{\Omega\cap\{\wuo > c_k\}} |\nabla\wuo|^{p-2} \nabla\wuo \nabla v_k =\int_\Omega |\nabla w_k|^{p}.
    \end{aligned}
\end{equation}
Combining \eqref{testfn},\eqref{sqnorm},\eqref{pnorm} we have that
\begin{equation*}
    \begin{aligned}
    \|w_k\|_{\wo}\leq \int_\Omega |\nabla\wuo|^{p-2} \nabla\wuo \nabla w_k &+ A_{q}(\wuo,w_k,\RR^{2N}).
    \end{aligned}
\end{equation*}
Note that, as stated in \cite{FP14}, for each $k\geq 1,$ the following holds:
\begin{equation*}
    \begin{aligned}
        \wuo <(2^k-1)w_{k-1}& \text{ in }\{w_k>0\},
 \\ &\text{ and }       
         \{w_k>0\}=\{\widetilde{\uo}>c_k\}\subseteq \{w_{k-1}>2^{k-1}\}.
    \end{aligned}
\end{equation*}
Together with the fact $w_k$ is a decreasing sequence, we get that 
\begin{equation}
    \begin{aligned}
     &\left[\int_\Omega C(\wuo^{-\delta}+ \wuo^{l})w_k\right] \leq   C\int_{\{w_{k}>0\}}[ c_{k}^{-\delta} w_{k-1} + (2^{k}-1)^{l-1}w_{k-1}^{l} ]
     \\
     & \leq C \int_{\{w_{k-1}>2^{-k}\}} \left[ c_{k}^{-\delta} 2^{k(l-1)}w_{k-1}^{l}+ (2^{k}-1)^{l-1}w_{k-1}^{l}\right]
     \\
     & \leq C \max \frac{1}{ 2^{kl}}\left\{ c_{k}^{-\delta} 2^{k(l-1)}, ({2^{k}-1})^{l-1}\right\}  2^{kl} \int_{\{w_{k-1}>2^{-k}\}} w_{k-1}^{l}.
    \end{aligned}
\end{equation}

Now we observe that for $\delta<1,$
\begin{equation*}
    \begin{aligned}
        \max \frac{1}{ 2^{kl}}\left\{ c_{k}^{-\delta} 2^{k(l-1)}, ({2^{k}-1})^{l-1}\right\} =\max \left\{\frac{(1-2^{-k})^{-\delta}}{2}, \frac{1}{2^{k}}(1-2^{-k})^{l-1}\right\}<1.
    \end{aligned}
\end{equation*}
So all together we obtain that
\begin{equation*}
    \begin{aligned}
        \int_{\Omega} |\nabla w_k|^{p} + \int \int_{\RR^{2N}}\frac{|w_{k}(x)-w_{k}(y)|^{q}}{|x-y|^{N+sq}}  < C 2^{kl}  \int_{\{w_{k-1}>2^{-k}\}} w_{k-1}^{l}.
    \end{aligned}
\end{equation*}
We also observe that 
\begin{equation*}
    \begin{aligned}
        \int_\Omega w_{k-1}^{l} \geq \int_{\{w_{k-1}>2^{-k}\}} w_{k-1}^{l} \geq 2^{-k l} |\{w_{k-1}>2^{-k}\}| \geq 2^{-kl}|\{w_{k}>0
        \}|.
    \end{aligned}
\end{equation*}
Now using the Sobolev and H\"older regularity, we get that
\begin{equation*}
    \begin{aligned}
        \|w_k\|_{L^{l}}^{l}&\leq \left(\int_{\Omega} w_{k}^{p^{*}}\right)^{\frac{l}{p^{*}}}|\{w_k>0\}|^{1-\frac{l}{p^{*}}}
     \leq S_{p}^{\frac{-l}{p}}\left(\int_\Omega|\nabla w_k|^{p}\right)^{\frac{l}{p}} |\{w_k>0\}|^{1-\frac{l}{p^{*}}}
        \\&< S_{p}^{\frac{-l}{p}} \left(C 2^{kl} \int w_{k-1}^{l}\right)^{\frac{l}{p}} 2^{(kl)(1-\frac{l}{p^{*}})} \left(\int w_{k-1}^{l}\right)^{1-\frac{qo}{p^{*}}}
        \\
        &
        \leq  S_{p}^{\frac{-l}{p}}  C^{\frac{l}{p}} 2^{kl(1+\frac{l}{p}-\frac{l}{p^{*}})} ( \|w_{k-1}\|_{L^{l}}^{l})^{(1+\frac{l}{p}-\frac{l}{p^{*}})}.
    \end{aligned}
\end{equation*}
Since $U_k:=\|w_k\|_{L^{l}}^{l}$ and $\frac{1}{p}-\frac{1}{p^{*}}=\frac{1}{n}$ we get that 
\begin{equation*}
    \begin{aligned}
        U_k \leq S_{p}^{\frac{-l}{p}}  C^{\frac{l}{p}} 2^{kl(1+\frac{l}{N})} U_{k-1}^{1+\frac{l}{N}}\leq C_0 2^{(k-1)(l+\frac{l^2}{N})} U_{k-1}^{1+\frac{l}{N}}
    \end{aligned}
\end{equation*} 
for chosen $C_0:=  S_{p}^{\frac{-l}{p}}  C^{\frac{l}{p}}   2^{l(1+\frac{l}{N})}.$
Moreover we chose $\alpha$ small enough independent of $k$ such that we get that 
\begin{equation*}
    \begin{aligned} U_0=\|\wuo\|_{L^{l}}^{l}=\alpha^{l}\|\uo\|_{L^{l}}^{l} \leq C_{0}^{\frac{-N}{l}} 2^{\frac{-N^2}{l^{2}}(l+\frac{l^2}{N})}.
    \end{aligned}
\end{equation*}
Then using \cite[Lemma 7.1]{Giu03} we get that
\begin{equation*}
    \begin{aligned}
        0=\lim_{k\rightarrow 0} U_k=\lim_{k\rightarrow 0} \int_\Omega (\wuo-c_k)^{l}_{+}=\int_\Omega (\wuo-1)^{l}_{+}.
    \end{aligned}
\end{equation*}
Hence, we get that $0\leq \uo\leq \frac{1}{\alpha}$ and conclude that $\uo\in L^{\infty}(\Omega).$
\end{proof}

Before proving the existence and uniqueness result for the sublinear perturbation problem, we will recall the Picone's identity for the both local and nonlocal cases.
\begin{theorem}\cite[Theorem 1.2]{BT20}\label{pcnloc}
Let $u>0$ and $u\geq 0$ be differentiable in $\Omega$. Then for  $1<r\leq t$ we have 
\begin{equation*}
    \begin{aligned}
       |\nabla u|^{t-2}\nabla u \nabla(\frac{v^{r}}{u^{r-1}})\leq \frac{r}{t} |\nabla v|^{t}+\frac{t-r}{t}|\nabla u|^{t}.
    \end{aligned}
\end{equation*}
If $r<t$, then the equality follows iff $u\equiv v.
$
\end{theorem}
\begin{theorem}\cite[Corollary 2.7]{GGM22}\label{pcnnonloc}
Let $0<s<1,$ $1<r<\infty$ and $1<t\leq r.$ Then for any $u,v>0$ two nonconstant measurable positive functions on $\Omega$ we have the following 
\begin{equation*}
    \begin{aligned}
       & [u(x)-u(y)]^{r-1}\left(\frac{u(x)^{t}-v(x)^{t}}{u(x)^{t-1}}-\frac{v(y)^{t}-u(y)^{t}}{v(y)^{t-1}}\right)
       \\
       &+  [v(x)-v(y)]^{r-1}\left(\frac{v(x)^{t}-u(x)^{t}}{v(x)^{t-1}}-\frac{u(y)^{t}-v(y)^{t}}{u(y)^{t-1}}\right) \geq 0
    \end{aligned}
\end{equation*}
holds for a.e. $x,y\in \Omega.$
\end{theorem}
Now, we give the proof of the main result of this section.
\\
\textit{\underline{Proof of Theorem \ref{extncuniqsub}}}
\\[3mm]
We can establish the existence of the solution in a manner analogous to that outlined in \cite[Theorem 2.7]{GKS21a}.

In the view of Lemma \ref{uobd} and the gradient H\"older regularity  result (Theorem \ref{gradhld}) we get that the solution $u \in C^{1,\alpha}(\Bar{\Omega})$ for some $\alpha \in (0,1).$ Only thing left to prove is the uniqueness of the solution. 
    \\
    Let $u_1$ and $u_2$ be two solutions of \eqref{uoeqn}. For $\epsilon>0,$ define $u_{i,\epsilon}=u_{i}+\epsilon, \thinspace i=1,2$ and choose 
\begin{equation*}
    \begin{aligned}
        \Phi:= \frac{u_{1,\epsilon}^{l+1}-u_{2,\epsilon}^{l+1}}{u_{1,\epsilon}^{l}}, \thickspace 
        \Psi:= \frac{u_{2,\epsilon}^{l+1}-u_{1,\epsilon}^{l+1}}{u_{2,\epsilon}^{l}}\in \wopo\cap\wsqo \text{ as test functions.}
    \end{aligned}
\end{equation*}
Clearly, 
\begin{equation}\label{pcneqn}
    \begin{aligned}
       & \int_{\Omega} |\nabla u_1|^{p-2}\nabla u_1 \nabla \Phi+ A_{q}(u_1,\Phi,\RR^{2N})=\int_{\Omega} (\frac{\lambda}{u_{1}^{\delta}}+u_{1}^{l})\Phi \text{ and }
       \\
        &\int_{\Omega} |\nabla u_2|^{p-2}\nabla u_2 \nabla \Psi+ A_{q}(u_2,\Psi,\RR^{2N})=\int_{\Omega} (\frac{\lambda}{u_{2}^{\delta}}+u_{2}^{l})\Psi.
    \end{aligned}
\end{equation}
Then, we conclude that 
\begin{equation}\label{pcnnoneg}
    \begin{aligned}
        (\frac{\lambda}{u_{1}^{\delta}}+u_{1}^{l})\Phi+ (\frac{\lambda}{u_{2}^{\delta}}+u_{2}^{l})\Psi 
     &=(u_{1,\epsilon}^{l+1}-u_{2,\epsilon}^{l+1})\left(\frac{\lambda}{u_{1}^{\delta} u_{1,\epsilon}^{l}}-\frac{\lambda}{u_{2}^{\delta} u_{2,\epsilon}^{l}}+\frac{u_{1}^{l}}{u_{1,\epsilon}^{l}}-\frac{u_{2}^{l}}{u_{2,\epsilon}^{l}}\right)
    \\ &\leq (u_{1,\epsilon}^{l+1}-u_{2,\epsilon}^{l+1})\left (\frac{u_{1}^{l}}{u_{1,\epsilon}^{l}}-\frac{u_{2}^{l}}{u_{2,\epsilon}^{l}}\right).
    \end{aligned}
\end{equation}
We get that $u_1, u_2 \in L^{\infty}(\Omega)$ using Lemma \ref{uobd}. Then, we can use the idea of \cite[pp.14]{GGM22} to pass through the limit on the right hand side of \eqref{pcneqn} as $\epsilon\rightarrow 0^{+}$.
 
Moreover, using the dominated convergence theorem and Fatou's lemma we get that 
\begin{equation}\label{pcncon}
    \begin{aligned}
               & \int_{\Omega} |\nabla u_{1,\epsilon}|^{p-2}\nabla u_{1,\epsilon} \nabla \frac{u_{1,\epsilon}^{l+1}-u_{2,\epsilon}^{l+1}}{u_{1,\epsilon}^{l}}+ A_{q}\left(u_{1,\epsilon},\frac{u_{1,\epsilon}^{l+1}-u_{2,\epsilon}^{l+1}}{u_{1,\epsilon}^{l}},\RR^{2N}\right)
               \\
               & 
               +  \int_{\Omega} |\nabla u_{2,\epsilon}|^{p-2}\nabla u_{2,\epsilon} \nabla \frac{u_{2,\epsilon}^{l+1}-u_{1,\epsilon}^{l+1}}{u_{2,\epsilon}^{l}}+ A_{q}\left(u_{2,\epsilon},\frac{u_{2,\epsilon}^{l+1}-u_{1,\epsilon}^{l+1}}{u_{2,\epsilon}^{l}},\RR^{2N}\right)
               \\
               &\xrightarrow{\epsilon\rightarrow 0^{+}}
               \\&  \int_{\Omega} |\nabla u_{1}|^{p-2}\nabla u_{1} \nabla \frac{u_{1}^{l+1}-u_{2}^{l+1}}{u_{1}^{l}}+ A_{q}\left(u_{1},\frac{u_{1}^{l+1}-u_{2}^{l+1}}{u_{1}^{l}},\RR^{2N}\right)
               \\
               & 
               +  \int_{\Omega} |\nabla u_{2}|^{p-2}\nabla u_{2} \nabla \frac{u_{2}^{l+1}-u_{1}^{l+1}}{u_{2}^{l}}+ A_{q}\left(u_{2},\frac{u_{2}^{l+1}-u_{1}^{l+1}}{u_{2}^{l}},\RR^{2N}\right).
    \end{aligned}
\end{equation}
Using Theorem \ref{pcnnonloc}, we obtain that 
\begin{equation} \label{pcnnonlc}
    \begin{aligned}
       A_{q}\left(u_{1},\frac{u_{1}^{l+1}-u_{2}^{l+1}}{u_{1}^{l}},\RR^{2N}\right)
              + A_{q}\left(u_{2},\frac{u_{2}^{l+1}-u_{1}^{l+1}}{u_{2}^{l}},\RR^{2N}\right) \geq 0.
    \end{aligned}
\end{equation}
Now, combining  \eqref{pcneqn},\eqref{pcnnoneg}, \eqref{pcncon}, and\eqref{pcnnonlc}  we get that 
\begin{equation*}
    \begin{aligned}
        \int_\Omega |\nabla u_1|^{p}+ |\nabla u_2|^{p} &\leq \int_{\Omega} |\nabla u_{1}|^{p-2}\nabla u_{1} \nabla \frac{u_{2}^{l+1}}{u_{1}^{l}} +  |\nabla u_{2}|^{p-2}\nabla u_{2} \nabla \frac{u_{1}^{l+1}}{u_{2}^{l}}.
    \end{aligned}
\end{equation*}
But using Theorem \ref{pcnloc} we have that
\begin{equation*}
    \begin{aligned}
        \int_{\Omega} |\nabla u_{1}|^{p-2}\nabla u_{1} \nabla \frac{u_{2}^{l+1}}{u_{1}^{l}} +  |\nabla u_{2}|^{p-2}\nabla u_{2} \nabla \frac{u_{1}^{l+1}}{u_{2}^{l}}
       \leq  \int_\Omega |\nabla u_1|^{p}+ |\nabla u_2|^{p}.
    \end{aligned}
\end{equation*}
Hence using Theorem \ref{pcnloc} we deduce that $u_1\equiv u_2.$
\hfill $\square$
 \section{Subcritical Superlinear Perturbation Problem}\label{subcrtical}
 In this section, we consider the superlinear problem 
 \begin{equation}\label{uoeqn1}
        \left\{
        \begin{array}{rll}
            \fp u +  \fqs u &= F(x,u);\;\; u>0 &\text{ in } \Omega,
         \\   u&=0 &\text{ in } \Omega^{c}
        \end{array}\right.
    \end{equation}
     where $F(x,s)= \lambda s^{-\delta}+s^{l-1}$  where $\max\{p,q\}\leq l< \min\{p^{*},q_{s}^{*}\},$ $0<\delta<1$ and $\lambda>0.$ 
     \subsection{Uniform Boundedness}
     We prove the uniform boundedness of the solution for the  subcritical perturbation taking inspiration from the proof of \cite[Theorem 2.3]{GKS20}.
 \begin{theorem}\label{uniformbd}
    Let $u\in\wo(\Omega)$ be a solution to the problem \eqref{uoeqn1}. Then $\|u\|_{L^{\infty}(\Omega)}$ is bounded depending on $\lambda, l,N,p$ and $\|u\|_{L^{l}(\Omega)}.$  
\end{theorem}
\begin{proof}
Without loss of generality, suppose $ u_{+} \neq 0 $. Let us choose $ \vartheta \geq \max\{1, \|u\|^{-1}_{L^{l}}\} $ and define $ v = (\vartheta \|u\|_{L^{l}(\Omega)})^{-1}u $. Then, $ v \in \wopo \cap \wsqo $ and $ \|v\|_{L^{l}} = \vartheta^{-1} $. For each $ k \in \mathbb{N} $, set $ w_k = (v - 1 + 2^{-k})_{+} $ and $ w_0 = v_{+} $. Clearly, $ w_k \in \wo $ and the sequence $ \{w_k\} $ is non-increasing, with $ w_k $ converging pointwise, as $k\to\infty$, to $ (v - 1)_{+} $ almost everywhere in $ \Omega $. Now, let $ U_k = \|w_k\|_{L^{l}}^{l} $, and observe that $ \{w_{k+1} > 0\} \subseteq \{w_k > 2^{-(k+1)}\} $, and $ v(x) < (2^{k+1} - 1)w_k(x) $ for $ x \in \{w_{k+1} > 0\} $. By Chebyshev's inequality, it follows that $ |\{w_{k+1} > 0\}| \leq 2^{(k+1)r} U_k $. Note that

    \begin{equation*}
        \begin{aligned}
            \|w_{k+1}\|_{\wsqo}^{q} =
\int_{\mathbb{R}^{2n}} \frac{|w_{k+1}(x) - w_{k+1}(y)|^q}{|x - y|^{N+sq}} \, dx \, dy &\leq A_{q}(v, w_{k+1},\RR^{2N})
\\
&=(\vartheta \|u\|_{L^{l}(\Omega)})^{1-q} A_{q}(u, w_{k+1},\RR^{2N}).
        \end{aligned}
    \end{equation*}
Observe that 
\begin{equation*}
    \begin{aligned}
        \int_{\Omega} |\nabla u|^{p-2} \nabla u \nabla w_{k+1}=\int_{\{v\geq 1-2^{-(k+1)}\} } (\vartheta \|u\|_{L^{l}(\Omega)})^{-1} |\nabla u|^{p} \geq 0.
    \end{aligned}
\end{equation*}
Hence,
\begin{equation}
    \begin{aligned}
       \|w_{k+1}\|_{\wsqo}^{q}    &\leq (\vartheta \|u\|_{L^{l}(\Omega)})^{1-q} A_{q}(u, w_{k+1},\RR^{2N})
       \\&\leq (\vartheta \|u\|_{L^{l}(\Omega)})^{1-q} \left( \int_{\Omega} |\nabla u|^{p-2} \nabla u \nabla w_{k+1}+A_{q}(u, w_{k+1},\RR^{2N})\right)
       \\&\leq C (\vartheta \|u\|_{L^{l}(\Omega)})^{1-q} \left[ \int_{\{w_{k+1}\geq 0\} } (u^{-\delta}+u^{l-1})w_{k+1}\right]
       \\&\leq C (\vartheta \|u\|_{L^{l}(\Omega)})^{1-q} \left[ \int_{\{w_{k+1}\geq 0\}\cap\{u\leq 1\} } u^{-\delta} w_{k+1}+ \int_{\{w_{k+1}\geq 0\}\cap\{u\geq 1\} }u^{l-1}w_{k+1}\right].
    \end{aligned}
\end{equation}
Noting that $(\vartheta \|u\|_{L^r(\Omega)})^{-1} \leq 1$, we have  $v(x) - 1 + 2^{-(k+1)} \leq 2^{-(k+1)}$ whenever $u(x) \leq 1,$
that is $w_{k+1}(x) \leq 2^{-(k+1)}$. Moreover from $\{w_{k+1}\geq 0\}$, we get that $u\geq \vartheta \|u\|_{L^{l}}(1-2^{-(k+1)}).$ Also, $u(x) = (\vartheta \|u\|_{L^{l}(\Omega)})v(x) < (\vartheta \|u\|_{L^{l}(\Omega)})(2^{k+1} - 1)w_k(x)$ in $\{w_{k+1}\geq 0\}.$
Thus, 
\begin{equation*}
    \begin{aligned}
        \|w_{k+1}\|_{\wsqo}^{q}    &\leq C (\vartheta \|u\|_{L^{l}(\Omega)})^{1-q}\Big[(\vartheta \|u\|_{L^{l}(\Omega)})^{-\delta}(1-2^{-(k+1)})^{-\delta}2^{-(k+1)}|\{w_{k+1}\geq 0\}
        \\
        &\hspace{1.0em}+ (\vartheta \|u\|_{L^{l}(\Omega)})^{l-1}(2^{k+1}-1)^{l-1}\int_{\{w_{k+1}\geq 0\}} w_{k}^{l}\Big]
        \\
        &\leq  C (\vartheta \|u\|_{L^{l}(\Omega)})^{1-q-\delta} (1-2^{-(k+1)})^{-\delta}2^{(l-1)(k+1)} U_{k}+C  (\vartheta \|u\|_{L^{l}(\Omega)})^{l-q} 2^{(l-1)(k+1)} U_{k}. 
    \end{aligned}
\end{equation*}
Observe that $1-2^{-(k+1)}>\frac{1}{2}$ and thus $((1-2^{-(k+1)})^{-\delta})\leq 2^{\delta} \leq 2$ and $\vartheta \|u\|_{L^{l}(\Omega)}\geq 1 $. Therefore, $\vartheta \|u\|_{L^{l}(\Omega)}^{1-q-\delta}\leq 1.$ Hence,
\begin{equation*}
    \begin{aligned}
        \|w_{k+1}\|_{\wsqo}^{q}    &\leq 2 C (\vartheta \|u\|_{L^{l}(\Omega)})^{l-q} 2^{(l-1)(k+1)} U_{k}.
    \end{aligned}
\end{equation*}
With the help of Sobolev embedding, we have
\begin{equation}
    \begin{aligned}
   U_{k+1} &= \|w_{k+1}\|_{L^l(\Omega)}^l \leq
\left(\int_{\Omega} w_{k+1}^{q^*_{s}} \, dx \right)^{l/q^{*}_{s}} |\Omega_{k+1}|^{1-l/q^{*}_{s}}
\\ & \leq (S^{-1} \|w_{k+1}\|_{\wsqo}^q )^{l/q} |\Omega_{k+1}|^{1 - l/q^*_{s}}
\\& \leq S^{-l/q} \left( 2C2^{(l-1)(k+1)}(\vartheta \|u\|_{L^l(\Omega)})^{l - q} U_k \right)^{l/q} ( 2^{(k+1)l}U_k )^{ \left( 1 - l/q^*_{s} \right)}
\\ & \leq C^k (\vartheta \|u\|_{L^l(\Omega)})^{\frac{l^{2}}{q
} - l} U_{k}^{1 + \frac{ls}{N}}
    \end{aligned}
\end{equation}
where $C > 1$ is independent of $k$. Let $\eta = C^{-N/(l s)} \in (0,1)$ and define $\gamma := \frac{l^{2} s}{N} + l - \frac{l}{q} > 0$. Choosing $\vartheta$ such that
\begin{equation*}
    \begin{aligned}
       \vartheta \geq \max \left\{ 1, \|u\|_{L^l(\Omega)}^{-1}, \left( \|u\|_{L^l(\Omega)}^{l^{2}/q - l} \eta^{-1}\right)^\frac{1}{\gamma}, C^\frac{ N^2}{(l s)^2} \right\},
    \end{aligned}
\end{equation*}
we claim that $ U_k \leq \frac{\eta^k}{\vartheta^l} $ and we will prove this by induction. For $ k = 0 $, we have
\begin{equation*}
U_0 = \|v_+\|_{L^l(\Omega)}^l \leq \|v\|_{L^l(\Omega)}^l = \vartheta^{-l}.
\end{equation*}
Now, assume that $ U_k \leq \frac{\eta^k}{\vartheta^l} $ holds for some $ k \in \mathbb{N} $. Then, we can write
\begin{equation*}
    \begin{aligned}
        U_{k+1} &\leq C^k (\vartheta \|u\|_{L^l(\Omega)})^{l^{2}/q - l} U^{1 + l s/N}_{k}
        \\& \leq C^k (\vartheta \|u\|_{L^l(\Omega)})^{l^{2}/q - l} \left( \frac{\eta^k}{\vartheta^l} \right)^{1 + l s/N}
        \\ & \leq \eta^k / \vartheta^l (\|u\|_{L^l(\Omega)})^{l^{2}/q - l} \vartheta^{-\gamma} \leq \eta^{k+1} / \vartheta^l.
    \end{aligned}
\end{equation*}
Therefore, the induction step holds, and the claim is proved. From this claim, it follows that $ U_k \to 0, $ as $ k \to \infty $, which implies $ w_k(x) \to 0 $, as $ k \to \infty $, almost everywhere in $ \Omega $. As a result, we have $ v(x) \leq 1 $ a.e. in $ \Omega $. Similarly, we can show that $ -v(x) \leq 1 $ a.e. in $ \Omega $. Therefore, $ \|v\|_{L^\infty(\Omega)} \leq 1 $ leading to $ \|u\|_{L^\infty(\Omega)} \leq \vartheta \|u\|_{L^l(\Omega)} $, which implies $ u \in L^\infty(\Omega) $.
\end{proof}
\subsection{Existence of Solution}
Let $\ullb$ be the solution of the following purely singular Dirichlet problem:
\begin{equation}\label{singulareqn}
\Bigg\{
    \begin{aligned}
       & \fp \ullb + \fqs \ullb =\lambda \ullb^{-\delta},\quad \ullb>0 \text{ in } \Omega,
        \\
        &\ullb=0\text{ in } \Omega^{c}.
        \end{aligned}
\end{equation}
The existence and the uniqueness of $\ullb$ are guaranteed by Corollary \ref{existenceofminimialsol} and Theorem \ref{wcp}. Using the Hopf's lemma \cite[Theorem 1.2]{AC23} we get the following lemma.
\begin{lemma}
    The problem \eqref{singulareqn} has a unique weak solution in $\wo(\Omega).$ Moreover $\ullb(x)\geq C d(x).$ 
\end{lemma}
Now we can modify the energy functional corresponding to \eqref{uoeqn1} using $\ullb$ to get the existence of one solution for \eqref{uoeqn1} for some $\lambda>0.$
\begin{lemma}
    There exists $\lambda_0>0$ such that \eqref{uoeqn1} with $\lambda=\lambda_0$ admits a solution  $u \in C^{1,\alpha}(\Bar{\Omega})$ for some $\alpha \in (0,1).$ 
\end{lemma}

\begin{sketch}
        Define
    \begin{equation*}
    \widetilde{f}(x,s)=\Bigg\{
        \begin{aligned}
            &\lambda s^{-\delta}+ s^{l} \text{ if } s\geq \ullb(x),\\
            & \lambda \ullb^{-\delta}+ \ullb^{l} \text{ if } s\leq \ullb(x).
        \end{aligned}
    \end{equation*}
    Let $\widetilde{F}(x,s)=\int_{0}^{s}  \widetilde{f}(x,t) dt,$ and $\Ellb:\wo(\Omega)\rightarrow \RR$ by 
    \begin{equation*}
        \begin{aligned}
\Ellb(u):=\frac{1}{p}\int_{\Omega} |\nabla u|^{p} + \frac{1}{q} [u]_{\wsqo}-\int_\Omega \widetilde{F}(x,u)\thinspace dx.
        \end{aligned}
    \end{equation*}
    Since the Strong Comparison Principle (Theorem \ref{scp}) and Theorem \ref{gradientholderregularity} are now established, we can modify  the calculations from \cite[Lemma 3.4]{GST07} within the framework of the mixed local-nonlocal operator to reach the desired conclusion.
\end{sketch}
Next we prove the existence of one solution for each $0<\lambda<\lambda_0.$
 \\
\textit{\underline{Proof of Theorem \ref{extsubcr}}}
\\[3mm]
Let $\ullb$ be the solution to \eqref{singulareqn} and $u_{\lambda_0}$ solves \eqref{uoeqn1} for $\lambda=\lambda_0.$ By Theorem \ref{gradientholderregularity}, both  $\ullb,u_{\lambda_0} \in C^{1,\alpha}(\Bar{\Omega}).$ Moreover, since the operator $\fp+\fqs$ is monotone, we deduce that $\ullb\leq u_{\lambda_0} $ adopting the calculations of \cite[pp.128]{GST07}. Set 
\begin{equation}
    \widetilde{f}_{\lambda}(x,s)=\left\{
    \begin{aligned}
        &\lambda u_{\lambda_0}^{-\delta}+u_{\lambda_0}^{l} \hspace{0.8em} \text{ if } s>u_{\lambda_0}(x),
        \\
        &\lambda s^{-\delta}+s^{l} \hspace{1.2em}\text{ if } \ullb(x)\leq s \leq u_{\lambda_0}(x),
        \\
        &
        \lambda \ullb^{-\delta}+\ullb^{l}\hspace{1.0em} \text{ if } s<\ullb(x).
    \end{aligned}
    \right.
\end{equation}
Let $\widetilde{F}_{\lambda}(x,s)=\int_{0}^{s} \widetilde{f}_{\lambda}(x,t) dt,$ and $\widetilde{\Ellb}:\wo(\Omega)\rightarrow \RR$ by 
    \begin{equation}
        \begin{aligned}
\widetilde{\Ellb}(u):=\frac{1}{p}\int_{\Omega} |\nabla u|^{p} + \frac{1}{q} [u]_{\wsqo}-\int_\Omega \widetilde{F}_{\lambda}(x,u)\thinspace dx.
        \end{aligned}
    \end{equation}
     Now $\widetilde{\Ellb}$ is coercive, weakly lower semi continuous and bounded below in $\wo(\Omega)$. Consequently $\widetilde{\Ellb}$ achieves its global minimum at some $\ulwt\in\wo(\Omega).$ Since $\widetilde{\Ellb}\in C^{1}(\wo(\Omega),\RR),$ it follows that $\ulwt$ solves the equation  $\fp \ulwt + \fqs \ulwt=\widetilde{f}_{\lambda}(x,\ulwt)$ in $\Omega.$ \\ Applying the strong maximum principle \cite[Theorem 1.2]{AC23} we get $\ulwt>0.$ From Theorem \ref{gradientholderregularity} we infer that $\ulwt\in C^{1,\alpha}(\overline{\Omega})$ and  Theorem \ref{scp}  further implies that $\ullb< \ulwt<u_{\lambda_0}.$ Therefore $\widetilde{f}_{\lambda}(x,\ulwt)=\lambda\ulwt^{-\delta}+\ulwt^{l},$ so  $\ulwt$ is a weak solution of \eqref{uoeqn1}.
\hfill $\square$\\
We conclude our paper with a final remark.
\begin{remark}
   Set $ \Lambda:=\sup\{\lambda>0: \eqref{uoeqn1}$ has a nontrivial solution $\}.$ It can be shown that whenever  $p=q,$ then $0<\Lambda<\infty.$  
   \end{remark}
Indeed, this remark can be reached by contradiction, using the isolation of the first eigenvalue $\lambda_1$ (See \cite[Theorem 1.1]{DFR19}). Suppose, there exists a sequence of $\lambda_{n}\rightarrow \infty$ such that \eqref{uoeqn1} admits a solution $u_n$ for each $n.$ There exists $\tilde{\lambda}_{0}>0$ such that the following holds:
   \begin{equation*}
       \begin{aligned}
           \lambda t^{-\delta} +t^{l}\geq (\lambda_1+\varepsilon) t^{p-1} , \text{ for all } t>0 \text{ and for any } 0<\varepsilon<1, \,\,  \lambda>\tilde{\lambda}_{0}.
       \end{aligned}
   \end{equation*}
For $\lambda_{n}>\tilde{\lambda}_{0},$ $u_n$ acts as a supersolution to the problem:
\begin{equation}
    \label{subsupeqn}
    \begin{aligned}
        \fp u + \fps u &= (\lambda_{1}+\varepsilon)u^{p-1}, \hspace{1.0em} u>0 \text{ in } \Omega, \text{ and } u|_{\Omega^{c}}=0,
    \end{aligned}
\end{equation}
for any $0<\varepsilon<1.$ Utilizing Hopf's lemma we can infer that $u_{n}\geq C(\tilde{\lambda}_{0}) \dis(x)$ independent of $n.$ Thus we can choose $\mu$ small enough such that $\mu<\lambda_{1}+\varepsilon$ and  $\mu \phi_{1}<C(\tilde{\lambda}_{0}) \dis(x)\leq u_{n}$ where $\phi_{1}$ is the (normalized) positive eigenfunction corresponding to $\lambda_1.$ Moreover $\mu \phi_{1}$ serves as a subsolution to \eqref{subsupeqn}. Applying the monotone iteration procedure yields a solution to \eqref{subsupeqn} for any $\varepsilon\in (0,1),$ leading to a contradiction since $\lambda_1$ is  an isolated point in the spectrum.

\bibliographystyle{plain}
	\bibliography{refnew}
\end{document}